\documentclass[letterpaper,11pt]{amsart}
\usepackage[margin=1.2in]{geometry}
\usepackage{amsmath,amsthm,amssymb}
\usepackage{xspace,xcolor}
\usepackage[breaklinks,colorlinks,citecolor=teal,linkcolor=teal,urlcolor=teal,pagebackref,hyperindex]{hyperref}
\usepackage[alphabetic]{amsrefs}
\usepackage[all]{xy}
\usepackage[english]{babel}
\usepackage{enumitem}
\usepackage{tikz,xcolor}
\usepackage{tikz-cd}
\usepackage{mathrsfs}
\usepackage{color}

\usepackage{mathtools}

\usepackage{version}
\excludeversion{hide}

\newcounter{intro}

\newtheorem{intro-conjecture}[intro]{Conjecture}
\newtheorem{intro-corollary}[intro]{Corollary}
\newtheorem{intro-theorem}[intro]{Theorem}

\newcommand{\theoremref}[1]{\hyperref[#1]{Theorem~\ref*{#1}}}
\newcommand{\lemmaref}[1]{\hyperref[#1]{Lemma~\ref*{#1}}}
\newcommand{\definitionref}[1]{\hyperref[#1]{Definition~\ref*{#1}}}
\newcommand{\propositionref}[1]{\hyperref[#1]{Proposition~\ref*{#1}}}
\newcommand{\conjectureref}[1]{\hyperref[#1]{Conjecture~\ref*{#1}}}
\newcommand{\corollaryref}[1]{\hyperref[#1]{Corollary~\ref*{#1}}}
\newcommand{\exampleref}[1]{\hyperref[#1]{Example~\ref*{#1}}}

\theoremstyle{plain}
\newtheorem{thm}{Theorem}[section]

\newtheorem{lem}[thm]{Lemma}
\newtheorem{prop}[thm]{Proposition}
\newtheorem{cor}[thm]{Corollary}

\theoremstyle{definition}
\newtheorem{defi}[thm]{Definition}

\newtheorem{eg}[thm]{Example}

\theoremstyle{remark}
\newtheorem{rmk}[thm]{Remark}

\def\Mustata{Mus\-ta\-\c{t}\u{a}\xspace}

\def\N{{\mathbf N}}
\def\Z{{\mathbf Z}}
\def\Q{{\mathbf Q}}
\def\R{{\mathbf R}}
\def\C{{\mathbf C}}

\def\A{{\mathbf A}}
\def\P{{\mathbf P}}

\def\cC{\mathcal{C}}
\def\cD{\mathcal{D}}
\def\cE{\mathcal{E}}

\def\cH{\mathcal{H}}

\def\cK{\mathcal{K}}

\def\cM{\mathcal{M}}

\def\cO{\mathcal{O}}

\def\cT{\mathcal{T}}

\def\cV{\mathcal{V}}

\def\.{\cdot}
\def\^{\widehat}

\def\de{\partial}

\def\({\left(}
\def\){\right)}

\renewcommand{\and}{ \ \ \text{ and } \ \ }

\DeclareMathOperator{\depth} {depth}

\begin{document}

\author[Q.~Chen]{Qianyu Chen}

\address{Department of Mathematics, University of Michigan, 530 Church Street, Ann Arbor, MI 48109, USA}

\email{qyc@umich.edu}

\author[B.~Dirks]{Bradley Dirks}

\address{Department of Mathematics, Stony Brook University, Stony Brook, NY 11794-3651, USA}

\email{bradley.dirks@stonybrook.edu}

\author[S.~Olano]{Sebasti\'{a}n Olano}

\address{Department of Mathematics, University of Toronto, 40 St. George St., Toronto, Ontario Canada, M5S 2E4}

\email{seolano@math.toronto.edu}

\thanks{Q.C. was partially supported by NSF Grant No. DMS-1952399, the Simons Collaboration grant Moduli of Varieties and AMS-Simons travel grant. B.D. was partially supported by the National Science Foundation under Grant No. DMS-1926686 and NSF-MSPRF grant DMS-2303070}

\title[Partial CCI varieties]{Partial Cohomologically Complete Intersections via Hodge Theory}

\begin{abstract} We introduce an invariant $c(X)$ associated to any complex algebraic variety $X$, which for varieties with isolated singularities measures the failure of dual Kodaira-Akizuki-Nakano vanishing. In general, it is characterized by Hodge--Lyubeznik numbers, the depth of Du Bois complexes, and the Hodge filtration on local cohomology modules.

We show that this invariant is computable in many examples, such as cones over rational homology manifolds or determinantal varieties. 
\end{abstract}

\maketitle

\section{Introduction}
Kodaira-Akizuki-Nakano vanishing for smooth, projective complex algebraic varieties is a cornerstone result used throughout complex algebraic geometry. It states that for a smooth projective complex variety $X$ and an ample line bundle $L$, we have the vanishing
\[ H^q(X, \Omega_X^p \otimes_{\cO_X} L) = 0 \text{ for all } p +q > \dim X,\]
and its dual version,
\[ H^q(X, \Omega_X^p \otimes_{\cO_X} L^{-1}) =0  \text{ for all } p+q < \dim X.\]

Though the general statement fails for singular varieties, \cite{GNAPGP}*{Thm. V.5.1} provides an analogue using Du Bois complexes $\underline{\Omega}_X^p$, which are Hodge-theoretic replacements for the K\"{a}hler differentials. Specifically, for a projective variety $X$ and an ample line bundle $L$ one has the hypercohomology vanishing
\[ \mathbb H^q(X, \underline{\Omega}_X^p \otimes_{\cO_X} L) = 0 \text{ for all } p +q > \dim X.\]
However, the dual version is known to fail even with the Du Bois complexes in place of the K\"{a}hler differentials.

In this paper, we introduce a Hodge-theoretic invariant $c(X) \in \Z_{\geq -1} \cup \{\infty\}$, defined by
\[ c(X) = \sup \{k\in \Z_{\geq -1} \mid {\rm depth}(\underline{\Omega}_X^p) \geq \dim X - p \text{ for all } p \leq k\},\]
which, in the isolated singularities case, is exactly the Hodge degree up to which the dual Kodaira-Akizuki-Nakano vanishing holds:
\begin{intro-theorem} \label{thm-Vanishing} Let $X$ be a projective variety with ample line bundle $L$. If $c(X) \geq k$, then
\[ \mathbb H^q(X, \underline{\Omega}_X^p \otimes_{\cO_X} L^{-1}) = 0 \text{ for all } p+q < \dim X, \, p \leq k.\]

If $X$ is equidimensional and $X_{\rm sing}$ is isolated, then the converse holds.
\end{intro-theorem}

Another version of this vanishing theorem, ranging over all $p$ but with a stronger bound on $p+q$, is \cite{PS}*{Thm. C}.

\begin{rmk} Note that if $c(X) \geq 0$, then $X$ is equidimensional. Suppose, for a contradiction, that $Z\subseteq X$ is an irreducible component with $\dim Z < \dim X$. Then for a closed point $x\in Z_{\rm reg}$ lying on no other irreducible component of $X$, we have
\[ {\rm depth}(\underline{\Omega}_X^0) \leq {\rm depth}_{\cO_{X,x}} (\underline{\Omega}_{X}^0)_x = \dim Z < \dim X.\]
\end{rmk}

\begin{rmk} \label{rmk-injThm} It is a consequence of the \emph{injectivity theorem} \cite{DBDeform} that
\[ {\rm depth}(\cO_X) \leq {\rm depth}(\underline{\Omega}_X^0),\]
and so, if $X$ is Cohen--Macaulay, then $c(X) \geq 0$.
\end{rmk}

As an application, \theoremref{thm-Vanishing} gives unobstructedness for deformations of certain mildly singular isolated Fano varieties, see \corollaryref{cor-Deformation} below.

Recall that the \emph{local cohomological defect} ${\rm lcdef}(X)$, defined by \cite{PS}, is a rough numerical measure of the variety's failure to be a local complete intersection. The condition ${\rm lcdef}(X) = 0$ is equivalent to Hellus and Schenzel's \emph{cohomologically complete intersection} (CCI) notion \cite{HellusSchenzel}. By \cite{MPLocCoh}*{Cor. 12.6}, we see that
\[ c(X) = + \infty \iff X \text{ is CCI},\]
hence the title of the paper.

Although $c(X)$ is defined in terms of the depth of the Du Bois complexes, and is motivated by a global vanishing theorem for projective $X$, it has three alternative characterizations related to local singularity invariants. The second condition of the theorem is written in terms of Saito's theory of mixed Hodge modules \cite{SaitoMHM}. The third characterization is through the Hodge--Lyubeznik numbers $\lambda^{p,q}_{r,s}(\cO_{X,x})$, defined by Garc\'{i}a-L\'{o}pez and Sabbah \cite{HodgeLyubeznik}.

\begin{intro-theorem} \label{thm-MaincX} Let $X$ be a complex algebraic variety. Then $c(X)$ is equal to
\begin{enumerate} 
\item \label{thm-LocCohCCI} $\sup\{k\mid F_k\cH^j_X(\cO_Y) = 0\text{ for all } j >{\rm codim}_Y(X)\},$ for $X\subseteq Y$ a closed embedding in a smooth variety $Y$.

\item \label{thm-depthDuBois}$\sup \{k \mid {\rm Gr}^F_{-p} {\rm DR}(\tau^{<0} \Q_X^H[\dim X]) = 0 \text{ for all } p \leq k\}$.

\item \label{thm-HLCCI} $\sup \{k \mid \lambda_{r,s}^{p,q}(\cO_{X,x}) = 0 \text{ for all } x\in X, s < \dim X, r, q\in \Z, p\geq -k\}$,
\end{enumerate}
\end{intro-theorem}

We call $\tau^{<0} \Q_X^H[\dim X] \in D^b({\rm MHM}(X))$ the \emph{CCI defect object}, and its support is precisely the closed subset where $X$ fails to be CCI, which we denote by $X_{\rm nCCI}$. 

\begin{hide}
Any rational homology manifold, meaning a variety such that the natural morphism $\Q_X[\dim X] \to {\rm IC}_X$ is a quasi-isomorphism, is CCI. Similarly, any connected LCI variety is CCI \cite{Dimca}*{Thm. 5.1.20}, hence the name.

Recall that the recently introduced invariant ${\rm HRH}(X)$ measures the Hodge-theoretic failure of the variety $X$ being a \emph{rational homology manifold}. We can also relate the partial versions of these singularity classes, as in the following theorem:
\begin{intro-corollary} \label{thm-HLHRH} Let $X$ be an equidimensional variety. Then
\[ {\rm HRH}(X) = \min\left\{c(X),\sup\{k \mid \lambda_{r,\dim X}^{p,q}(\cO_{X,x}) = {\rm I}\lambda_r^{p,q}(\cO_{X,x}) \text{ for all } x\in X, q,r\in \Z, p\geq -k\}\right\}.\]
\end{intro-corollary}

In particular, we have the inequality
\[ {\rm HRH}(X) \leq c(X).\]

\begin{eg} This inequality can be strict: any hypersurface which is not a rational homology manifold gives such an example. For example, $f = xy-zw$ satisfies ${\rm HRH}(X) = 0$ but $c(X) = +\infty$.
\end{eg}
\end{hide}

We have the following theorem which is analogous to the corresponding inequality for the non-rational homology manifold locus \cite{DOR1}*{Thm. C}. In \cite{DOR1}, a ``generic variant'' ${\rm lcdef}_{\rm gen}(X)$ of the ${\rm lcdef}(X)$ invariant was defined. In this paper, we define another generic variant ${\rm lcdef}_{\rm gen}^{>0}(X)$ which only concerns the local cohomology sheaves $\cH^{c+j}_X(\cO_Y)$ with $j > 0$. See Subsection \ref{ss-Inequality} below for the precise definition.

\begin{intro-theorem} \label{thm-inequalityCCI} Assume $\underline{\Omega}_X^0$ is Cohen--Macaulay, in other words, $c(X) \geq 0$. Then we have the inequality
\[c(X) \leq \frac{{\rm codim}_X(X_{\rm nCCI}) - {\rm lcdef}_{\rm gen}^{>0}(X) - 3}{2}.\]
\end{intro-theorem}

As a simple consequence, using the fact that ${\rm lcdef}_{\rm gen}^{>0}(X) >0$ if and only if $X_{\rm nCCI} \neq \emptyset$, we see that if $\underline{\Omega}_X^0$ is Cohen--Macaulay, then
\[ {\rm codim}_X(X_{\rm nCCI}) \geq 2c(X) + 4.\]

\begin{rmk} In low dimensions, this has an interesting consequence. If $\dim X \leq 3$, then $\underline{\Omega}_X^0$ being Cohen--Macaulay implies $\Q_X[\dim X]$ is perverse.

In particular, using Remark \ref{rmk-injThm}, we can recover a consequence of a theorem of Dao--Takagi \cite{DaoTakagi}*{Cor. 2.8} (pointed out to us by Hyunsuk Kim and Mihnea Popa, to whom we are thankful). If $X$ is Cohen--Macaulay and $\dim X \leq 3$, then $\Q_X[\dim X]$ is perverse.
\end{rmk}

We stress that $c(X)$ is computable in many classes of examples: for a CCI variety $X$ (for example, a connected LCI variety or a variety with quotient singularities), we have $c(X) = +\infty$. The inequality of \theoremref{thm-inequalityCCI} will be used in Section \ref{sec-Determinantal} below to study its value for determinantal varieties.

The following theorem provides a further class of examples. The key takeaway is that, for varieties with a one-step resolution, many of the relevant singularity invariants are computable in terms of the Hodge structure on primitive cohomology of the exceptional divisor. This is a consequence of the more general \theoremref{thm-general}. The general version is used in \cite{CDORSecant} to understand the Hodge-theoretic singularity invariants of secant varieties and in particular, to compute $c(X)$ and ${\rm HRH}(X)$ in that setting. The statement is given in terms of Saito's mixed Hodge modules, though we remark that it applies to theories of mixed sheaves as in Section \ref{subsect-MixedSheaves} below.

Let $X$ be an irreducible complex algebraic variety of dimension at least $2$. Let $f \colon \widetilde{X} \to X$ be a projective morphism which is an isomorphism over $X\setminus \{x\}$, for some point $x\in X$ and such that $\widetilde{X}$ and $\widetilde{Z} = f^{-1}(x)$ are irreducible rational homology manifolds, with $d = \dim \widetilde{Z} = \dim \widetilde{X} - 1$. Assume moreover that $\iota \colon X \hookrightarrow Y$ is a closed embedding into a smooth variety $Y$ with $c= \dim Y - \dim X$, and write $i\colon \{x\} \hookrightarrow Y$.

\begin{intro-theorem} \label{thm-IntroGeneral} In the above notation, we have
\begin{enumerate} 
    \item For all $0 < j \leq \dim X -2$, we have an isomorphism of pure Hodge modules of weight $\dim Y + c +j+1$:
    \[ \cH^{c+j}_{X}(\cO_Y)  \cong \mathbf{V}^{d-j}.\]
    \item If ${\rm Gr}^W_w \cH^c_{X}(\cO_Y) \neq 0$, then $w \in \{\dim Y + c,\dim Y+c+1\}$, and we have an isomorphism
    \[ {\rm Gr}^W_{\dim Y + c+1} \cH^c_{X}(\cO_Y) \cong \mathbf{V}^d.\]
 \item ${\rm lcdef}(X) \leq \dim X -2$.
\end{enumerate}
Here $\mathbf{V}^{d-j}= i_*H_{\rm prim}^{d-j}(\widetilde{Z})(-c-j-1).$
\end{intro-theorem}

The proof is quite formal: for weight reasons, in this set-up, many connecting maps in the pull-back long exact sequence associated to the map $\Q_X \to f_* \Q_{\widetilde{X}}$ must vanish. Then we use the Hard Lefschetz theorem on $\widetilde{Z}$ to conclude.

\begin{rmk} These formulas generalize and refine Hartshorne--Polini's \cite{HartshornePolini} computation of the local cohomology modules of cones over smooth projective varieties. See Section \ref{sec-Contractions} for more details.
\end{rmk}

\begin{hide}
\begin{intro-theorem} \label{thm-IntroGeneral} In the above notation, we have
\begin{enumerate} \item ${\rm lcdef}(X) \leq \dim X -2$.
    \item For all $0 < j \leq \dim X -2$, we have an isomorphism of pure Hodge modules of weight $\dim X +j+1$:
    \[ \cH^{c+j}_{X}(\cO_Y) \cong \iota_*\cH^j \mathbf D^H_X(-c) \cong i_*H_{\rm prim}^{d-j}(\widetilde{Z})(-c-j-1).\]
    \item We have ${\rm Gr}^W_w \cH^0\mathbf D_X^H \neq 0$ implies $w \in \{\dim X,\dim X+1\}$, and an isomorphism
    \[ {\rm Gr}^W_{\dim Y + c+1} \cH^c_{X}(\cO_Y) \cong \iota_*({\rm Gr}^W_{\dim X+1} \cH^0 \mathbf D_X^H)(-c) \cong i_* H^{d}_{\rm prim}(\widetilde{Z})(-c-1).\]

\item We have an isomorphism in $D^b ({\rm MHM}(Y))$:
\[ \iota_* f_* \Q_{\widetilde{X}}^H[\dim X] \cong \iota_*{\rm IC}_X^H \oplus  i_* H^{\dim X}(\widetilde{Z}) \oplus \bigoplus_{\ell =1}^{d-1} \left(i_*H^{\dim X+\ell}(\widetilde{Z})[-\ell] \oplus i_*H^{\dim X+\ell}(\widetilde{Z})(\ell)[\ell]\right).\]
\end{enumerate}
\end{intro-theorem}
\end{hide}

\theoremref{thm-IntroGeneral} gives the following control over singularity invariants:

\begin{intro-corollary} \label{cor-CollectResults} In the above set-up, the following can be written in terms of primitive Hodge numbers of $\widetilde{Z}$:
\begin{enumerate} \item(\corollaryref{cor-kCCI}) $c(X)$,
\item (\corollaryref{cor-HRHIsolated}) ${\rm HRH}(X)$,
\item (\corollaryref{cor-GenLevel}) For a closed embedding $X\subseteq Y$ of codimension $q$ into a smooth variety, the generating level of the Hodge filtration on ${\rm Gr}^W_{\dim X+2q+1}\cH^q_X(\cO_Y)$ and that on $\cH^{q+j}_X(\cO_Y)$ for $j>0$,
\item (\theoremref{thm-HLIsolated}) $\lambda_{r,s}^{p,q}(\cO_{X,x})$,
\item (\theoremref{thm-HLIsolated}) ${\rm I}\lambda_{r}^{p,q}(\cO_{X,x})$.
\end{enumerate}
\end{intro-corollary}

\medskip

The paper is organized by first developing the invariant $c(X)$ and its basic properties, then studying it in determinantal examples, and finally computing it and related invariants for varieties admitting one-step resolutions, in particular, for contractions of zero sections.

\noindent {\bf Acknowledgments} The authors would like to thank Bhargav Bhatt, James Hotchkiss, Hyunsuk Kim, Mircea \Mustata, Sung Gi Park, Mihnea Popa, Debaditya Raychaudhury, Claude Sabbah, Christian Schnell, Rosie Shen, Sridhar Venkatesh, Duc Vo and Jakub Witaszek for many conversations about ideas occurring in this work. 

\section{The Partial CCI Invariant} \label{sec-PartialCCI}
Throughout the paper, $X$ is a connected complex algebraic variety.

We will use the following standard notation from the theory of mixed Hodge modules: $\Q_X^H$ is the constant object on $X$ and $\mathbf D_X^H = \mathbf D_X(\Q_X^H[\dim X])(-\dim X)$ is the ``dual'' constant object, both lying in $D^b({\rm MHM}(X))$. Here $(-\dim X)$ is the Tate twist functor, which, in terms of underlying bifiltered $\cD$-modules, simply shifts the filtrations.

We will also make use of the functor ${\rm Gr}^F_p {\rm DR}(-) \colon D^b({\rm MHM}(X)) \to D^b_{\rm coh}(\cO_X)$. For a single filtered right $\cD$-module $(\cM,F)$ on a smooth variety $Y$, we have
\[ {\rm Gr}^F_p {\rm DR}(\cM) = {\rm Gr}^F_{p-\dim Y}(\cM) \otimes \bigwedge^{\dim Y} \cT_Y \xrightarrow[]{\nabla} \dots \xrightarrow[]{\nabla} {\rm Gr}^F_p(\cM).\]

 \subsection{Equivalent Characterizations} \label{subsec-EquivalentCharacterizations}
 Recall from the introduction that a variety $X$ is called a \emph{cohomologically complete intersection} (CCI) if $\Q_X[\dim X]$ is a perverse sheaf. We define the \emph{CCI defect object of $X$} to be
 \[ \tau^{<0} \Q_X^H[\dim X] \in D^b({\rm MHM}(X)).\]

It is convenient to use this object to define $c(X)$, and then to prove that it agrees with the characterization in terms of the depth of Du Bois complexes given in the introduction.

\begin{defi} We define
\[ c(X) = \sup\{k \in \Z_{\geq -1} \mid {\rm Gr}^F_{-p} {\rm DR}(\cH^j \Q_X^H[\dim X]) = 0 \text{ for all } j < 0, p \leq k\}\]
or, equivalently, by duality,
\[ c(X) = \sup \{k \in \Z_{\geq -1} \mid {\rm Gr}^F_{p-\dim X} {\rm DR}(\cH^j \mathbf D_X^H) = 0 \text{ for all } j > 0, p\leq k\}.\]
\end{defi}
If $X \subseteq Y$ is embedded into a smooth variety $Y$, then this condition can be characterized in terms of local cohomology (compare \cite{DOR1}*{Thm. D(1)}). First, we can rewrite local cohomology in terms of $\mathbf D_X^H$.
\begin{lem}\label{lem-DXLocCoh} Let $i\colon X \to Y$ be a closed embedding with $\dim Y - \dim X = c$, where $Y$ is smooth. Then there is a natural isomorphism of bi-filtered $\cD_Y$-modules underlying mixed Hodge modules
\[ i_* \cH^{j} \mathbf D_X^H (-c) \cong (\cH^{c+j}_X(\cO_Y),F,W),\]
where on the right, the object is the local cohomology $\cD_Y$-module of the structure sheaf $\cO_Y$ along $X$.
\end{lem}
\begin{proof} By definition, $(\cH^q_X(\cO_Y),F,W)$ underlies the mixed Hodge module $\cH^q i_* i^! \Q_Y^H[\dim Y]$. As $Y$ is smooth, it is a rational homology manifold, so this is isomorphic to $\cH^q i_* i^! \mathbf D_Y^H$. Using the fact that $i^* \Q_Y^H \cong \Q_X^H$, it is easy to see that
\[ \cH^q( i_* \mathbf D^H_X (-c)[-c]) = i_* \cH^{q-c} \mathbf D_X^H (-c),\]
as claimed.
\end{proof}

More generally, if $f \colon X_1 \to X_2$ is any morphism of complex algebraic varieties, then
\begin{equation} \label{eq-DXRule} f^! \mathbf D_{X_2}^H \cong \mathbf D_{X_1}^H(\dim X_1 -\dim X_2)[\dim X_1 - \dim X_2].\end{equation}

The following proposition proves the equivalence between \theoremref{thm-MaincX}\eqref{thm-LocCohCCI} and \theoremref{thm-MaincX}\eqref{thm-depthDuBois}.

\begin{prop} Let $i \colon X \to Y$ be an embedding of the complex algebraic variety $X$ into a smooth variety $Y$. Let $c = \dim Y - \dim X$. Then $c(X)\geq k$ if and only if $F_k \cH^{c+j}_X(\cO_Y) = 0$ for all $j > 0$.
\end{prop}
\begin{proof} By \lemmaref{lem-DXLocCoh} above, we have
\[ i_* {\rm Gr}^F_{p-\dim X} {\rm DR}_X(\cH^j(\mathbf D_X^H)) = {\rm Gr}^F_{p-\dim Y} {\rm DR}_Y(\cH_X^{c+j}(\cO_Y)),\]
where the difference in index is due to the Tate twist. So the claim is clear.
\end{proof}

We now turn to the equivalence between this definition of the $c$-invariant and that from the introduction. Recall that for an object $K \in D^b_{\rm coh}(\cO_X)$, we define the \emph{depth} of $K$ at a closed point $x\in X$ to be
\[ {\rm depth}_x(K) = \min \{i \mid \cH^{-i} {\rm RHom}_{\cO_{X,x}}(K_x,\omega_{X,x}^\bullet) \neq 0\},\]
and we let
\[ \depth(K) = \min_{x\in {\rm supp}(K)} \depth_x(K)\]
which, if $K$ is concentrated in degree $0$, recovers the usual notion of depth. 

So we have
\[ {\rm depth}(K) \geq j \text{ if and only if } {\rm RHom}_{\cO_{X,x}}(K_x,\omega_{X,x}^\bullet) \in D^{\leq -j}_{\rm coh}(\cO_{X,x}) \text{ for all } x\in {\rm supp}(K),\]
which gives the equivalence
\[ {\rm depth}(K) \geq j \text{ if and only if } \mathbb D(K) \in D^{\leq -j}_{\rm coh}(\cO_X).\]

We will use the mixed Hodge module interpretation of the Du Bois complex \cite{SaitoMHC} (see also \cite{MOPW}*{Lem. 2.3}), which says that for all $p \in \Z$, we have a quasi-isomorphism
\[ {\rm Gr}^F_{-p} {\rm DR}(\Q_X^H) \cong \underline{\Omega}_X^p[-p].\]

\begin{prop} For any complex algebraic variety $X$, we have $c(X)\geq k$ if and only if
\[ {\rm depth}(\underline{\Omega}_X^p) \geq \dim X - p \text{ for all } p\leq k.\]
\end{prop}
\begin{proof} For any $p$, we have the spectral sequence
\[ {}^p E_2^{i,j} = \cH^i {\rm Gr}^F_p {\rm DR}(\cH^j \mathbf D(\Q_X^H[\dim X])) \implies \cH^{i+j} {\rm Gr}^F_p {\rm DR}(\mathbf D(\Q_X^H[\dim X])).\]

Note that 
\[{\rm Gr}^F_p {\rm DR}(\mathbf D(\Q_X^H[\dim X])) = R \cH om_{\cO_X}({\rm Gr}^F_{-p} {\rm DR}(\Q_X^H[\dim X]),\omega_X^\bullet)\]
\[ = R \cH om_{\cO_X}(\underline{\Omega}_X^p, \omega_X^\bullet)[p-\dim X].\]

Thus, ${\rm depth}(\underline{\Omega}_X^p) \geq \dim X - p$ is equivalent to the vanishing 
\[ \cH^\ell {\rm Gr}^F_p{\rm DR}(\mathbf D(\Q_X^H[\dim X])) =0 \text{ for all } \ell  > 0 .\]

If $c(X)\geq k$, then for all $p\leq k$, we have ${}^p E_2^{i,j} = 0$ for all $j > 0$, and ${}^p E_2^{i,\ell} = 0$ for all $i > 0, \ell \in \Z$ (which is a general property for the de Rham complex of a mixed Hodge module). Thus, for any $i+j > 0$, we get ${}^p E_2^{i,j} = 0$ and therefore, $ {}^p E_\infty^{i,j} = 0$, which implies $\cH^{i+j} {\rm Gr}^F_p {\rm DR}(\mathbf D(\Q_X^H[\dim X])) = 0$. As noted above, this implies ${\rm depth}(\underline{\Omega}_X^{p})\geq \dim X -p$.

We prove the converse by induction on $k$. So assume $c(X)\geq k-1$ and that ${\rm depth}(\underline{\Omega}_X^k) \geq \dim X -k$. By the assumption on $c(X)$, we know that ${}^k E^{i,j}_2 = 0$ for $j > 0$ and $i\neq 0$. Indeed, by assumption, $k$ is the first possibly non-zero index of the Hodge filtration on $\cH^j \mathbf D(\Q_X^H[\dim X])$, and for such an index, it is known that ${\rm Gr}^F_k {\rm DR}(\cH^j\mathbf D(\Q_X^H[\dim X]))$ is a coherent sheaf in degree $0$.

Thus, under the assumption that $c(X)\geq k-1$, the spectral sequence degenerates at $E_2$. We get, then, that ${}^k E_2^{0,j} = {}^k E_\infty^{0,j}$ is a subquotient of $\cH^j {\rm Gr}^F_k {\rm DR}(\mathbf D(\Q_X^H[\dim X]))$. By the equivalence mentioned above, ${\rm depth}(\underline{\Omega}_X^k) \geq \dim X -k$ implies that this is $0$ for $j > 0$. As this was the only possibly non-zero cohomology of ${\rm Gr}^F_k {\rm DR}(\cH^j(\mathbf D(\Q_X^H[\dim X])))$, this proves that the complex is quasi-isomorphic to zero, hence $c(X) \geq k$.
\end{proof}

This depth condition has gained significant interest recently in the work of \cites{MPLocCoh, PS}. Burke \cite{Burke} proves more refined relations between the Hodge filtration on local cohomology and depth of Du Bois complexes.

\subsection{Dual KAN vanishing and deformations} \label{subsec-DualKAN}
Having identified $c(X)$ both with the Hodge filtration on local cohomology and with depth conditions on the Du Bois complexes, we now turn to the global vanishing theorem that originally motivates the invariant.

Note that the following theorem is a refinement of \theoremref{thm-Vanishing}. Indeed, the forward direction is identical. For the converse as stated in \theoremref{thm-Vanishing}, we observe that if $X$ is equidimensional and $X_{\rm sing}$ is isolated, then
\[ X_{\rm nCCI} \subseteq X_{\rm sing}\]
is also a finite subset of $X$.

\begin{thm} Let $X$ be a projective variety with ample line bundle $L$. If $c(X) \geq k$, then
\[ \mathbb H^q(X, \underline{\Omega}_X^p \otimes_{\cO_X} L^{-1}) = 0 \text{ for all } p+q < \dim X, \, p \leq k.\]

If $X_{\rm nCCI}$ is isolated, then the converse holds.
\end{thm}
\begin{proof} Use some power of the ample line bundle $L$ to embed $X$ into projective space $i \colon X \to \P^N$, and write $c = N - \dim X$. By abuse of notation, we will not use $i_*(-)$ for complexes on $\P^N$ that are supported on $X$.

For any $p \leq c(X)$, we have an isomorphism
\[ {\rm Gr}^F_{p-N} {\rm DR}_{\P^N}((\cH^{c} i_* i^! \Q_{\P^N}^H[N])) \cong {\rm Gr}^F_{p-N} {\rm DR}_{\P^N}(i_* i^! \Q_{\P^N}^H[N])[c].\]

Hence, we conclude

\[ {\rm Gr}^F_{p-N} {\rm DR}_{\P^N}((\cH^c i_* i^! \Q_{\P^N}^H[N])) \cong \mathbb D( \underline{\Omega}_X^p)[p-\dim X].\]

But $\cH^c i_* i^! \Q_{\P^N}^H[N] \in {\rm MHM}(X)$, so we conclude the desired vanishing by Saito's vanishing theorem \cite{SaitoMHM}*{Prop. 2.33}, see also \cite{SchnellVanishing}*{Thm. 2.1}.

For the converse, assume $X_{\rm nCCI} = S \subseteq X$ is a finite subset with $\iota \colon S \to \P^N$ the closed embedding. In particular, we have an exact triangle
\[ \cH^c i_* i^! \Q_{\P^N}^H[N] \to i_* i^! \Q_{\P^N}^H[N+c] \to \iota_* \cC \xrightarrow[]{+1},\]
where $\cC \in D^b{\rm MHM}(S) = D^b{\rm MHS}_\Q$. By definition, $\iota_* \cC = \tau^{\geq 1}i_* i^! \Q_{\P^N}^H[N+c]  \in D^{\geq 1}({\rm MHM}(\P^N))$. 

It is easy to see that $c(X) \geq k$ is equivalent to the vanishing
\[ {\rm Gr}^F_{p-N} {\rm DR}(\tau^{\geq 1}i_* i^! \Q_{\P^N}^H[N+c]) = 0 \text{ for all } p \leq k,\]
or equivalently,
\[ \iota_* {\rm Gr}^F_{p-N}\cC = 0 \text{ for all } p \leq k,\]
and so we focus on proving this vanishing.

If we apply $L \otimes_{\cO} {\rm Gr}^F_{p-N} {\rm DR}(-)$ to the above exact triangle, we get the exact triangle
\[ L \otimes_{\cO} {\rm Gr}^F_{p-N} {\rm DR}(\cH^c i_* i^! \Q_{\P^N}^H[N]) \to L \otimes_{\cO} \mathbb D(\underline{\Omega}_X^{p})[p-\dim X] \to \iota_* {\rm Gr}^F_{p-N} \cC \xrightarrow[]{+1}\]
in $D^b_{\rm coh}(\cO_{\P^N})$. By Saito's vanishing, the left object satisfies
\[ H^q(L\otimes_{\cO} {\rm Gr}^F_{p-N} {\rm DR}(\cH^c i_* i^! \Q_{\P^N}^H[N])) = 0 \text{ for all } q > 0,\]
and so the long exact sequence in hypercohomology gives isomorphisms
\[ H^{p+q-\dim X}(L \otimes_{\cO} \mathbb D(\underline{\Omega}_X^p)) \cong \cH^q {\rm Gr}^F_{p-N} \cC \text{ for all } p \leq k, q > 0.\]

We now see immediately that the vanishing statement is equivalent to $c(X) \geq k$, as desired.
\end{proof}

This vanishing theorem can be used to prove unobstructedness of deformations for certain mildly singular Fano varieties. The proof is essentially the same as \cite{FL1}*{Thm. 4.5}. See also \cite{Tenie}*{Thm. 3.4} for a related result.

Recall that the Du Bois complexes, by definition, come with canonical morphisms
\[ \Omega_X^p \to \underline{\Omega}_X^p\]
in $D^b_{\rm coh}(\cO_X)$, for all $p\in \Z$.

\begin{cor} \label{cor-Deformation} Let $X$ be a Gorenstein Fano variety with isolated singularities, so $X$ is projective and $\omega_X$ is an anti-ample line bundle.

If $\dim X \geq 4$ and $c(X) \geq 1$, then $\cH^{\dim X - 3} \underline{\Omega}_X^1 = 0$ if and only if ${\rm Ext}^2(L_{X/\C}, \cO_X) = 0$. In this case, deformations of $X$ are unobstructed. 

If $\dim X = 3$, then the natural morphism $\Omega_X^1\to \cH^0 \underline{\Omega}_X^1$ is surjective if and only if ${\rm Ext}^2(L_{X/\C}, \cO_X) = 0$. In this case, deformations of $X$ are unobstructed. 
\end{cor}
\begin{proof} Let $X$ be a $d$-dimensional Fano variety with Gorenstein singularities. If $d \geq 4$, assume moreover $c(X) \geq 1$. Note that if $\dim X = 3$ and $X$ is Cohen--Macaulay, then Remark \ref{rmk-injThm} gives $c(X) \geq 0$, and then by \theoremref{thm-inequalityCCI}, we conclude that $X$ is CCI, so $c(X) \geq 1$ in this case, too.

Consider the exact triangle
\[ L_{X/\C} \to \underline{\Omega}_X^1 \to \cC^\bullet \xrightarrow[]{+1}\]
defined via the natural maps through K\"{a}hler differentials. 

We want to study the vanishing ${\rm Ext}^2( L_{X/\C}, \cO_X) = 0$. By Grothendieck--Serre duality, it suffices to prove $H^{d-2}(L_{X/\C} \otimes \omega_X) = 0$. Twisting the exact triangle by $\omega_X$ and applying hypercohomology, we get
\[ R\Gamma(L_{X/\C} \otimes \omega_X) \to R\Gamma(\underline{\Omega}_X^1 \otimes \omega_X) \to R\Gamma(\cC^\bullet \otimes \omega_X) \xrightarrow[]{+1}.\]

We have an exact sequence
\[ \mathbb H^{d-3}(\underline{\Omega}_X^1 \otimes \omega_X) \to  \mathbb H^{d-3}(\cC^\bullet \otimes \omega_X) \to H^{d-2}(L_{X/\C} \otimes \omega_X) \to \mathbb H^{d-2}(\underline{\Omega}_X^1 \otimes \omega_X) = 0,\]
the rightmost term vanishing by dual Akizuki--Nakano vanishing \theoremref{thm-Vanishing}. Hence, we see that our desired vanishing is implied by the vanishing
\[ \mathbb H^{d-3}(\cC^\bullet \otimes \omega_X) = 0.\]

We have $\mathbb H^{d-3} (\underline{\Omega}_X^1 \otimes \omega_X) = 0$ by \theoremref{thm-Vanishing}, and so the desired vanishing is \emph{equivalent} to the vanishing $\mathbb H^{d-3}(\cC^\bullet \otimes \omega_X)=0$.

If we assume moreover that $\cC^\bullet$ is supported on a finite subset (for example, if $X$ has isolated singularities), then this vanishing is equivalent to the vanishing $\cH^{d-3} \cC^\bullet = 0$. As $L_{X/\C} \in D^{\leq 0}_{\rm coh}(\cO_X)$, we have, if $d\geq 4$, the isomorphism
\[ \cH^{d-3} \underline{\Omega}_X^1 \cong \cH^{d-3} \cC^\bullet\]
and if $d =3$, the exact sequence
\[ \Omega_X^1 \to \cH^0 \underline{\Omega}_X^1 \to \cH^0 \cC^\bullet \to 0,\]
from which the result is clear.
\end{proof}

\subsection{Relation with ${\rm HRH}(X)$} \label{subsec-HRHRelation}
We now return to some general properties of the invariant $c(X)$, specifically how it relates to the recently introduced invariant ${\rm HRH}(X)$. As the CCI property is a natural weakening of the rational homology manifold condition, such a comparison is expected.

An equidimensional variety $X$ admits a natural exact triangle in $D^b({\rm MHM}(X))$
\[  \cK_X^\bullet \to \Q_X^H[\dim X]\to {\rm IC}_X^H \xrightarrow[]{+1},\]
where the leftmost object is called the \emph{RHM defect object}, because its support
\[ {\rm Supp}(\cK_X^\bullet) = X_{\rm nRS}\]
is the smallest closed subset of $X$ whose complement is a rational homology manifold.

Recall the following definition \cites{PPLefschetz,DOR1}:
\begin{defi} Let $X$ be an equidimensional complex algebraic variety. Then
\[ {\rm HRH}(X) = \sup\{k \mid {\rm Gr}^F_{-p} {\rm DR}(\cK_X^\bullet) =0 \text{ for all } p \leq k\}.\]
\end{defi}

The condition ${\rm HRH}(X) \geq m$ is called the $D_m$ property in \cites{PPLefschetz,PPFactorial}.

This definition closely mirrors that of $c(X)$ given in \theoremref{thm-MaincX}\eqref{thm-depthDuBois}. There is also a description in terms of local cohomology in the spirit of \theoremref{thm-MaincX}\eqref{thm-LocCohCCI}, see \cite{DOR1}*{Thm. D(1)}. We immediately get the following:

\begin{lem} \label{lem-diffCCIHRH} \cite{DOR1}*{Thm. D(1)} We have ${\rm HRH}(X) \leq c(X)$. Moreover, equality holds if and only if ${\rm Gr}^F_{-p} {\rm DR}(\cH^0 \cK_X^\bullet) = 0$ for all $p\leq c(X)$.
\end{lem}

Now we turn to a more topological invariant. The \emph{local cohomological defect} of a variety $X$ was introduced in \cite{PS} and related to the depth of the Du Bois complexes in \cite{MPLocCoh}. Here, we define this defect to be
\[ {\rm lcdef}(X) = \max\{j \mid \cH^{-j}\Q_X^H[\dim X] \neq 0\},\]
and so we see ${\rm lcdef}(X) = 0$ if and only if $X$ is CCI.

We now prove \theoremref{thm-MaincX}\eqref{thm-HLCCI}. First, we give an equivalent definition of Hodge--Lyubeznik numbers to that in \cite{HodgeLyubeznik}. Let $x\in X$ be a point in a possibly singular variety $X$. We define
\[ \lambda_{r,s}^{p,q}(\cO_{X,x}) = \dim_{\C} {\rm Gr}^F_{- p} {\rm Gr}^W_{p+q} \cH^r_x(\cH^{\dim X -s} \mathbf D_X^H(\dim X)),\]
where $\cH^r_x = \cH^r i_x^!$ and $i_x\colon \{x\} \to X$ is the inclusion of the point.

It is not hard to see that this agrees with the definition in \emph{loc. cit.} by \lemmaref{lem-DXLocCoh} above.

We also define the \emph{intersection Hodge--Lyubeznik numbers}
\[ {\rm I}\lambda_{r}^{p,q}(\cO_{X,x}) = \dim_{\C} {\rm Gr}^F_{- p} {\rm Gr}^W_{p+q} \cH^r_x({\rm IC}_X^H(\dim X)),\]
which do not depend on $s$, as ${\rm IC}_X^H$ has a unique non-vanishing cohomology module.

We first observe that these numbers are trivial if we look at a point where $X$ is CCI. Recall from the introduction the notation $X_{\rm nCCI}$ which denotes the non-CCI locus of $X$.

\begin{lem} \label{lem-HLRHM} Let $x\in X\setminus X_{\rm nCCI}$. Then
\[ \lambda_{\dim X,\dim X}^{0,0}(\cO_{X,x}) = 1\]
and all other Hodge--Lyubeznik numbers vanish.
\end{lem}
\begin{proof} Note that $i_x^!$ factors through the open embedding $U = X\setminus X_{\rm nCCI} \to X$, so we can assume from the beginning that $X$ is CCI.

Then $\mathbf D_X^H(\dim X) = \cH^0 \mathbf D_X^H(\dim X)$ by definition of CCI. This shows that $\lambda_{r,s}^{p,q}(\cO_{X,x}) = 0$ for $s \neq \dim X$. For $s =\dim X$, we have
\[ i_x^! \cH^0 \mathbf D_X^H(\dim X) = i_x^! \mathbf D_X^H(\dim X) = \mathbf D_{\{x\}}^H[-\dim X],\]
where we have used the relation \eqref{eq-DXRule}. But then $\mathbf D_{\{x\}}^H \cong \Q^H$ is the trivial Hodge structure, which proves the claim.
\end{proof}

Recall that for $M^\bullet \in D^b({\rm MHM}(X))$, we define
\[ p(M^\bullet) = \min\{p \mid {\rm Gr}^F_p {\rm DR}(M^\bullet) \neq 0\},\]
and it is not hard to see that we have equality
\begin{equation} \label{eq-pOfComplex} p(M^\bullet) = \min_{i \in \Z} p(\cH^i M^\bullet). \end{equation}

If $(V,F,W)$ is a mixed Hodge structure, then it is clear that
\begin{equation} \label{eq-pmixedGrW} p(V) = \min_{w\in \Z} p({\rm Gr}^W_w V),\end{equation}
as $W$ is an exhaustive filtration and the Hodge filtration on ${\rm Gr}^W_w V$ is induced by the filtration on $V$.

The proof of \theoremref{thm-MaincX}\eqref{thm-HLCCI} follows quickly using the following result \cite{DOR1}*{Lem. 1.8}.

\begin{lem} \label{lem-pRestrict} Let $M^\bullet \in D^b({\rm MHM}(X))$. Then
\[ p(M^\bullet) = \min_{x\in X} p(i_x^! M^\bullet),\]
where $i_x \colon \{x\} \to X$ is the inclusion of the point.
\end{lem}

\begin{proof}[Proof of \theoremref{thm-MaincX}\eqref{thm-HLCCI}] We have by definition
\[ c(X) \geq k \iff {\rm Gr}^F_{p-\dim X} {\rm DR}(\cH^j \mathbf D_X^H) = 0 \text{ for all } p \leq k, j >0,\]
which is true if and only if
\[ p(\cH^j \mathbf D_X^H(\dim X)) \geq k+1 \text{ for all } j > 0.\]

By \lemmaref{lem-pRestrict} above, this is equivalent to
\[ \text{ for all } x\in X_{\rm nCCI}, \, p(i_x^!\cH^j \mathbf D_X^H(\dim X)) \geq k+1 \text{ for all } j > 0\]

By the equality \eqref{eq-pOfComplex} and some reindexing, this is equivalent to
\[ \text{ for all } x\in X_{\rm nCCI}, \, p(\cH^{r}i_x^!\cH^{\dim X-s} \mathbf D_X^H(\dim X)) \geq k+1 \text{ for all } r \in \Z, s < \dim X\]

Now, $\cH^r i_x^! \cH^{\dim X -s} \mathbf D_X^H(\dim X)$ is a mixed Hodge structure, and by applying \eqref{eq-pmixedGrW} we see that\small
\[ p(\cH^r i_x^! \cH^{\dim X-s} \mathbf D_X^H(\dim X)) \geq k+1 \iff \text{ for all }w\in \Z, p({\rm Gr}^W_{w} \cH^r i_x^!\cH^{\dim X-s}\mathbf D_X^H(\dim X)) \geq k+1.\]
\normalsize

We have shown that $c(X) \geq k$ if and only if for all $x\in X_{\rm nCCI}, s < \dim X, r\in \Z$ and $q \in \Z$, we have
\[ {\rm Gr}^F_j {\rm Gr}^W_{p+q} \cH^r i_x^! \cH^{\dim X-s} \mathbf D_X^H(\dim X) =0 \text{ for all } j \leq k,\]
or equivalently,
\[ {\rm Gr}^F_{-p} {\rm Gr}^W_{p+q} \cH^r i_x^! \cH^{\dim X-s} \mathbf D_X^H(\dim X) =0 \text{ for all } p \geq -k,\]
which is equivalent to the claimed vanishing of Hodge--Lyubeznik numbers.
\end{proof}

We also obtain the following result, characterizing ${\rm HRH}(X)$ in terms of the intersection Hodge--Lyubeznik numbers and the invariant $c(X)$.

\begin{cor} \label{thm-HLHRH} For any equidimensional complex algebraic variety $X$, we have equality
\[ {\rm HRH}(X) = \min\left\{c(X),\sup \Lambda_X\right\},\]
where
\[\Lambda_X = \{k \in \Z \mid \text{ for all } x\in X, {\rm I}\lambda_r^{p,q}(\cO_{X,x}) = \begin{cases} 0 & \text{ for all } p \geq -k, q,r \in \Z, (p,q,r) \neq (0,0,\dim X) \\ 1 & (p,q,r) = (0,0,\dim X) \end{cases}\}.\]
\end{cor}

\begin{proof} For $k \in \Z_{\geq -1}$ fixed, we prove that ${\rm HRH}(X) \geq k$ if and only if the right hand side of the claimed equality is $\geq k$. In particular, in either case, $c(X) \geq k$, so we can assume this throughout the proof.

By \lemmaref{lem-diffCCIHRH}, ${\rm HRH}(X) \geq k$ if and only if
\[ {\rm Gr}^F_{p-\dim X} {\rm DR}(\cH^0\mathbf D_X^H/ {\rm IC}_X^H) = 0 \text{ for all } p \leq k,\]
or by Tate twisting,
\[ {\rm Gr}^F_{p} {\rm DR}((\cH^0\mathbf D_X^H/ {\rm IC}_X^H)(\dim X)) = 0 \text{ for all } p \leq k.\]

Thus, ${\rm HRH}(X) \geq k$ is equivalent to
\[ p((\cH^0\mathbf D_X^H/ {\rm IC}_X^H)(\dim X)) \geq k+1,\]
which by \lemmaref{lem-pRestrict} is equivalent to the inequality for all $x\in X_{\rm nRS}$:
\[p(i_x^!(\cH^0\mathbf D_X^H/ {\rm IC}_X^H)(\dim X)) \geq k+1.\]

By applying the exact functor $F_k(-)$ to the object $i_x^!(\cH^0\mathbf D_X^H/ {\rm IC}_X^H)(\dim X) \in D^b({\rm MHS})$, this inequality is equivalent to the vanishing for all $x\in X_{\rm nRS}$
\[ F_k(i_x^!(\cH^0\mathbf D_X^H/ {\rm IC}_X^H)(\dim X)) = 0.\]

By the exact triangle in $D^b({\rm MHS})$
\[ i_x^! {\rm IC}_X^H(\dim X) \to i_x^! \cH^0 \mathbf D_X^H (\dim X) \to i_x^!(\cH^0\mathbf D_X^H/ {\rm IC}_X^H)(\dim X) \xrightarrow[]{+1},\]
if we apply the exact functor $F_k(-)$ to the long exact sequence in cohomology we conclude that ${\rm HRH}(X) \geq k$ if and only if for all $x\in X_{\rm nRS}$, the natural map
\begin{equation} \label{eq-qisCondition} F_k(i_x^! {\rm IC}_X^H(\dim X)) \to F_k(i_x^! \cH^0 \mathbf D_X^H (\dim X)) \text{ is a quasi-isomorphism}.\end{equation}

Using the assumption that $c(X) \geq k$, we see that the right hand side is equal to $F_k i_x^! \mathbf D_X^H(\dim X) = F_k \Q^H[-\dim X]$. Thus, the condition \eqref{eq-qisCondition} is equivalent to
\[\cH^j F_k i_x^! {\rm IC}_X^H(\dim X) = 0 \text{ for all }j < \dim X,\]
\[\cH^{\dim X} F_k i_x^! {\rm IC}_X^H(\dim X) \cong F_k \Q^H.\]

In terms of associated graded pieces, this is equivalent to
\[ {\rm Gr}^F_{-p} {\rm Gr}^W_{p+q} \cH^ji_x^! {\rm IC}_X^H(\dim X) = 0 \text{ for all }j < \dim X, p \geq -k\]
\[{\rm Gr}^F_{-p} {\rm Gr}^W_{p+q}\cH^{\dim X} i_x^! {\rm IC}_X^H(\dim X) = 0, 0 > p \geq -k\]
\[{\rm Gr}^F_0 {\rm Gr}^W_{0} \cH^{\dim X}i_x^! {\rm IC}_X^H(\dim X) \cong \C.\]

For the last claim, note that, indeed, the $1$-dimensionality is equivalent to the map
\[\cH^{\dim X} F_k i_x^! {\rm IC}_X^H(\dim X) \to F_k \Q^H\]
being an isomorphism, because this map is always surjective. It is dual to the natural injective morphism
\[ \Q^H \to \cH^{-\dim X} i_x^* {\rm IC}_X^H.\]
\end{proof}

We collect here some useful facts about equidimensional complex algebraic varieties satisfying ${\rm HRH}(X) \geq k$. First, we define the \emph{generic local cohomological defect}, following \cite{DOR1}. We define $d(i) = \dim {\rm supp} \cH^{-i} \cK_X^\bullet$, so that $\dim X_{\rm nRS} = \max_{i\geq 0} d(i)$.

Then we have
\[ {\rm lcdef}_{\rm gen}(X) = \begin{cases} 0 & X \text{ is an RHM}  \\ \max\{i \mid d(i) = \dim(X_{\rm nRS})\} & X \text{ is not an RHM}\end{cases}.\]

\begin{prop} Let $X$ be an algebraic variety.

\begin{enumerate} \item $X$ is a rational homology manifold if and only if ${\rm HRH}(X) \geq k$ for all $k$, in which case we write ${\rm HRH}(X) = +\infty$.

\item \cite{PPLefschetz}*{Cor. 7.5} If $X$ satisfies ${\rm HRH}(X) \geq k$ for some $k \geq 0$, then $X$ is a rational homology manifold away from a closed subset of codimension at least $2k+3$.

\item If $X$ has Du Bois singularities, then it has rational singularities if and only if ${\rm HRH}(X) \geq 0$.

\item \cite{DOR1}*{Thm. C} If ${\rm HRH}(X) \geq 0$, then we have an inequality
\[ {\rm lcdef}_{\rm gen}(X) + 2 {\rm HRH}(X) + 3 \leq {\rm codim}_X(X_{\rm nRS}),\]
where $X_{\rm nRS} = {\rm supp}(\cK_X^\bullet)$ is the locus of closed points where $X$ is not a rational homology manifold.
\end{enumerate}
\end{prop}

\subsection{Non-CCI Locus and Codimension Bounds} \label{ss-Inequality}
The comparison with ${\rm HRH}(X)$ suggests that an analogous inequality could hold with $c(X)$ in place of ${\rm HRH}(X)$, and an adjusted ${\rm lcdef}_{\rm gen}(X)$ definition, as in \theoremref{thm-inequalityCCI}.

The following is immediate from the definition.

\begin{lem} We have
\[ c(X) = \inf\{p \in \Z_{\geq -1} \mid {\rm Gr}^F_{p+1-\dim X} {\rm DR}(\tau^{>0}\mathbf D_X^H) \neq 0\} \]
and
\[ X_{\rm nCCI} = {\rm supp}(\tau^{>0}\mathbf D_X^H) = {\rm supp}(\tau^{<0} \Q_X^H[\dim X]).\]
\end{lem}

In analogy with the above, we define
\[ {\rm lcdef}^{>0}_{\rm gen}(X) = \begin{cases} 0 & X\text{ is CCI} \\ \max (\{0\} \cup\{i> 0 \mid d(i) = \dim(X_{\rm nCCI})\}) & X \text{ is not CCI}\end{cases}.\]

The following is easily seen:

\begin{lem} If $X$ is not CCI, then the condition $\dim X_{\rm nCCI} =\dim X_{\rm nRS}$ is equivalent to ${\rm lcdef}_{\rm gen}(X) >0$, in which case we have ${\rm lcdef}_{\rm gen}(X) = {\rm lcdef}_{\rm gen}^{>0}(X)$.
\end{lem}

As a result, if the equality in the lemma holds, then the inequality of \theoremref{thm-inequalityCCI} is stronger than that of \cite{DOR1}*{Thm. C}, due to the inequality $c(X) \geq {\rm HRH}(X)$.

We now explain the proof of \theoremref{thm-inequalityCCI}. The basic support dimension bound is given in the following lemma, related to \cite{SaitoOnTheHodge}*{Pg. 10 Rmk.(ii)}. The right hand side of the following inequality is called the \emph{de Rham level}.

\begin{lem} \label{lem-LevelDimSupp} Let $M$ be a non-zero mixed Hodge module on a variety $X$. Let $p(M)$ be defined as above, and let $\rho(M) = \max\{p \mid {\rm Gr}^F_{p} {\rm DR}(M) \neq 0\}$.

Then
\[ \dim {\rm Supp}(M) \leq \rho(M) - p(M).\]
\end{lem}
\begin{proof} There exists a smooth, connected open subset $U\subseteq {\rm Supp}(M)$ with $\dim U = \dim {\rm Supp}(M)$ and such that $M\vert_U$ is a mixed Hodge module associated to an admissible variation of mixed Hodge structure. As the ${\rm Gr}^F_p {\rm DR}(-)$ commutes with restriction to open subsets, we have the inequalities
\[ p(M) \leq p(M\vert_U) \leq \rho(M\vert_U) \leq \rho(M),\]
thus,
\[ \rho(M\vert_U) - p(M\vert_U) \leq \rho(M) - p(M).\]

So we can restrict to the case that $X$ is smooth and $(\cV,F)$ is the filtered coherent $\cO_U$-module underlying $M$. Let $F_{p(\cV)} \cV$ be the first non-zero Hodge piece and assume $q \geq p(\cV)$ is the largest $q$ such that ${\rm Gr}^F_q(\cV) \neq 0$, which exists since the Hodge filtration on a variation of mixed Hodge structures is finite.

Then ${\rm Gr}^F_{p(\cV)} {\rm DR}(\cV) = {\rm Gr}^F_{p(\cV)} \cV \neq 0$ and ${\rm Gr}^F_{q+\dim U} {\rm DR}(\cV) = {\rm Gr}^F_q(\cV) \otimes_{\cO_U} \bigwedge^{\dim U} \cT_U [\dim U]$ are both non-zero. Thus,
\[ \rho(\cV) - p(\cV) = q+\dim U - p(\cV) \geq \dim U,\]
as claimed.
\end{proof}

Recall that \cite{PPLefschetz}*{Prop. 6.4} showed the following about the RHM defect object for a purely $n$-dimensional variety $X$:
\begin{equation} \label{eq-RHMDefect} \cK_X^\bullet \text{ has weights }\leq n -1.\end{equation}

Hence, for any $j > 0$, we have
\[ \cH^{-j} \cK_X^\bullet = \cH^{-j}\Q_X^H[\dim X] \text{ has weights } \leq n-j-1.\]

We get the following dimension bound:
\begin{lem} \label{lem-DimSuppBoundcInvariant} Let $X$ be purely $n$-dimensional. For any $j> 0$, we have
\[ \dim {\rm Supp}(\cH^{-j}\Q_X^H[\dim X]) \leq n - j - 2c(X) - 3.\]
\end{lem}
\begin{proof} As noted above, $\cH^{-j}\Q_X^H[\dim X]$ has weights $\leq n-j-1$. By definition of the $c$-invariant, we have
\[ {\rm Gr}^F_{-p} {\rm DR}( \cH^{-j}\Q_X^H[\dim X]) = 0 \text{ for all } p \leq c(X),\]
in other words, $\rho(\cH^{-j}\Q_X^H[\dim X]) < -c(X)$.

Combining the weight bound and this vanishing, \cite{PPLefschetz}*{Prop. 5.2} gives the vanishing
\[ {\rm Gr}^F_{-p} {\rm DR}(\cH^{-j}\Q_X^H[\dim X]) = 0 \text{ for all } p \geq n - j -1 -c(X),\]
in other words, $p(\cH^{-j} \Q_X^H[\dim X]) \geq c(X)+j+2-n$. Thus, we get (by \lemmaref{lem-LevelDimSupp}) the inequality
\[ \dim {\rm Supp}(\cH^{-j}\Q_X^H[\dim X]) \leq (-c(X) -1) - (c(X) +j+2-n) = n - j - 2c(X) - 3.\]
\end{proof}

We now prove \theoremref{thm-inequalityCCI}.

\begin{proof}[Proof of \theoremref{thm-inequalityCCI}] The condition $\underline{\Omega}_X^0$ is Cohen--Macaulay is equivalent to $c(X) \geq 0$. Hence, the proof is immediate from \lemmaref{lem-DimSuppBoundcInvariant}: choose $j>0$ maximal so that the value of $\dim {\rm supp} \cH^{-j} \Q_X^H[n]$ is maximal. This maximal dimension agrees with $\dim X_{\rm nCCI}$ and this $j$ is, by definition, ${\rm lcdef}_{\rm gen}^{>0}(X)$. Rearranging the inequality, we get
\[ {\rm lcdef}_{\rm gen}^{>0}(X) + 2c(X) + 3 \leq \dim X - \dim X_{\rm nCCI} = {\rm codim}_X(X_{\rm nCCI}),\]
as desired.
\end{proof}

\begin{cor} If $c(X) \geq 0$, then $X$ is CCI away from a subset of codimension $\geq 2c(X) + 4$.

In particular, if $\dim X \leq 3$, then $\underline{\Omega}_X^0$ being Cohen--Macaulay implies $\Q_X[\dim X]$ is perverse.
\end{cor}
\begin{proof} The first claim immediately implies the second, keeping in mind $\underline{\Omega}_X^0$ being Cohen--Macaulay is equivalent to $c(X) \geq 0$ by \theoremref{thm-MaincX}\eqref{thm-depthDuBois}.

For the first claim, either $X$ is CCI in which case the claim is obvious. Otherwise, if $X$ is not CCI, then ${\rm lcdef}_{\rm gen}^{>0}(X) \geq 1$, so the claim follows from \theoremref{thm-inequalityCCI}.
\end{proof}

\section{Examples: Determinantal Varieties}
\label{sec-Determinantal}
In this section, we show how the main theorem of \cite{RW} (and the discussion concerning it in \cite{DOR1}*{Sec. 10}) can be used to compute ${\rm lcdef}_{\rm gen}^{>0}(-)$ and to bound $c(Z)$ when $Z$ is a determinantal variety. This also illustrates the utility of the general inequality \theoremref{thm-inequalityCCI} when one can compute the topological quantities ${\rm codim}_X(X_{\rm nCCI})$ and ${\rm lcdef}_{\rm gen}^{>0}(X)$.

We will be interested in subvarieties defined by rank conditions in the following spaces: \begin{enumerate} \item (Generic) $X = {\rm Mat}_{m,n}(\C)$ with $m\geq n$, \item (Odd skew) $X = {\rm Mat}_{n}(\C)^{\rm skew}$, $n$ odd, \item  (Even skew) $X = {\rm Mat}_{n}(\C)^{\rm skew}$, $n$ even, \item (Symmetric) $X = {\rm Mat}_{n}(\C)^{\rm sym}$. \end{enumerate}

Following \cite{DOR1}*{Sec. 10}, in cases (1) and (4), we let $Z_p$ denote the subvariety of matrices of rank $\leq p$ and in cases (2) and (3), we let $Z_p$ denote the subvariety of matrices of rank $\leq 2p$. 

Following \cite{RW}, we let $D_p$ be the intersection homology $\cD_X$-module associated to the trivial local system on $Z_{p,\rm reg}$. We let $\Gamma(X)$ denote the Grothendieck group of holonomic $\cD_X$-modules. For $p$ fixed, we write
\[ H_p(q) = \sum_{j\geq 0} \left[\cH_{Z_p}^j(\cO_X)\right]\cdot q^j \in \Gamma(X)[q].\]

Finally, for $a\geq b\geq 0$, let $\binom{a}{b}_q$ be the $q$-binomial coefficient, defined by
\[ \binom{a}{b}_q = \frac{(1-q^a)\dots (1-q^{a-b+1})}{(1-q^b)\dots (1-q)}.\]

We state here the main result of \cite{RW}, giving a formula for the polynomial $H_p(q) \in \Gamma(X)[q]$.

\begin{thm}[\cite{RW}*{Main Thm.}]\label{thm-RW} In the notation above, we have the following formula for $H_p(q)$ in the cases (1)-(4).
\begin{enumerate} \item(Generic) For all $0\leq p < n$, we have \[H_p(q) = \sum_{s=0}^p [D_s]\cdot q^{(n-p)^2+(n-s)(m-n)} \binom{n-s-1}{p-s}_{q^2}.\]
\item(Odd skew) Write $n= 2m+1$. Then for all $0\leq p< m$, we have \[H_p(q) = \sum_{s=0}^p [D_s]\cdot q^{2(m-p)^2+(m-p) +2(p-s)} \binom{m-1-s}{p-s}_{q^4}.\]
\item(Even skew) Write $n = 2m$. Then for all $0\leq p < m$, we have \[H_p(q) = \sum_{s=0}^p [D_s]\cdot q^{2(m-p)^2-(m-p)} \binom{m-1-s}{p-s}_{q^4}.\]
\item(Symmetric) For all $0\leq p < n$, we have \[H_p(q) = \sum_{\ell = 0}^{\lfloor \frac{p}{2}\rfloor} [D_{p-2\ell}]\cdot q^{1+\binom{n-p+2\ell+1}{2} - \binom{2\ell+2}{2}}\binom{\lfloor \frac{n-p+2\ell-1}{2}\rfloor}{\ell}_{q^{-4}}.\]
\end{enumerate}
\end{thm}

As $Z_0$ is a point (hence smooth), we assume in the first three cases that $p\geq 1$. In case (4), it is known that $Z_1$ is a rational homology manifold, hence CCI, and so in case (4) we assume $p \geq 2$.

Recall that if ${\rm lcdef}_{\rm gen}(Z_p) > 0$, then we have equality 
\[ {\rm lcdef}_{\rm gen}(Z_p) = {\rm lcdef}_{\rm gen}^{>0}(Z_p).\]

By \cite{DOR1}*{Prop. 10.6}, we have the following:
\begin{enumerate} \item(Generic) ${\rm lcdef}_{\rm gen}(Z_p) = m+n-2p-2$.

\item(Odd skew) ${\rm lcdef}_{\rm gen}(Z_p) = 4(m-p-1)+2$.

\item(Even skew) ${\rm lcdef}_{\rm gen}(Z_p) = 4(m-p-1)$.

\item(Symmetric) ${\rm lcdef}_{\rm gen}(Z_p) = 2(n-p-1)$ (we assume $p\geq 2$).
\end{enumerate}

This immediately leads to the following:

\begin{cor} \label{cor-computelcdefgenpos} In the notation above, we have
\begin{enumerate} \item(Generic) ${\rm lcdef}_{\rm gen}^{>0}(Z_p) = m+n-2p-2$. In particular, for $m = n = p+1$, we know that $Z_p$ has hypersurface singularities, hence is CCI.

\item(Odd skew) ${\rm lcdef}_{\rm gen}^{>0}(Z_p) = 4(m-p-1)+2$.

\item(Even skew) ${\rm lcdef}_{\rm gen}^{>0}(Z_p) = 4(m-p-1)$. In particular, for $p = m-1$, we know that $Z_p$ has hypersurface singularities, hence is CCI.

\item(Symmetric) We have ${\rm lcdef}_{\rm gen}^{>0}(Z_p) = 2(n-p-1)$ (we assume $p\geq 2$). In particular, for $p = n-1$, we know $Z_p$ has hypersurface singularities, hence is CCI.
\end{enumerate}
\end{cor}
\begin{proof} It suffices to study the cases when ${\rm lcdef}_{\rm gen}(Z_p) = 0$ in the above formula. This equality is impossible in Case (2).

For Case (1), that equality is only possible for $m= n = p+1$, but in that case $Z_p$ is a hypersurface in $X$, hence CCI.

For Case (3), that equality is possible only for $m = p+1$, but in that case, $Z_p$ is a hypersurface in $X$ by \cite{DOR1}*{Cor. 10.3}, hence CCI.

For Case (4), again, equality is only possible for $n = p+1$, but in that case $Z_p$ is a hypersurface in $X$ by \cite{DOR1}*{Cor. 10.3}.
\end{proof}

\begin{cor} \label{cor-computeZnCCI} We have
\begin{enumerate} \item(Generic) $Z_{p,{\rm nCCI}} = \begin{cases} \emptyset & m=n=p+1 \\ Z_{p-1} & \text{otherwise}\end{cases}$,

\item(Odd skew) $Z_{p,{\rm nCCI}} = Z_{p-1}$,

\item(Even skew) $Z_{p,{\rm nCCI}} = \begin{cases} \emptyset & m=p+1 \\ Z_{p-1} & \text{otherwise}\end{cases}$,

\item(Symmetric) $Z_{p,{\rm nCCI}} = \begin{cases} \emptyset & n=p+1 \\ Z_{p-2} & \text{otherwise}\end{cases}$ (we assume $p\geq 2$).
\end{enumerate}
\end{cor}

Finally, this allows us to obtain a result similar to that of \cite{DOR1}*{Prop. 10.10}, using the inequality \theoremref{thm-inequalityCCI}. Note that since all varieties considered are Cohen--Macaulay (because they have rational singularities), we know $c(Z_p) \geq 0$.

\begin{cor} In the notation above, we have
\begin{enumerate} \item(Generic) For $m=n=p+1$, we have $c(Z_p) = \infty$. Otherwise, $c(Z_p) = 0$.

\item(Odd skew) $c(Z_p) \in \{0,1\}$.

\item(Even skew) For $m = p+1$, $c(Z_p) = \infty$. Otherwise, $c(Z_p) \in \{0,1\}$.

\item(Symmetric) For $n = p+1$, $c(Z_p) = \infty$. Otherwise, $c(Z_p) \in \{0,1\}$ (when $p\geq 2$).
\end{enumerate}
\end{cor}

\section{One-step Resolutions and Primitive Cohomology}\label{sec-MainThm}
In this section, we provide a technical result that allows for the computation of $c(X)$ and related singularity invariants when the variety $X$ admits a one-step resolution to a rational homology manifold. Moreover, as our argument essentially only uses six functors and the weight formalism for mixed Hodge modules, it generalizes to the following setting.

\subsection{$A$-mixed sheaves} \label{subsect-MixedSheaves}
Let $A\subseteq \R$ be a subfield. Throughout this section, we work in a theory of $A$-mixed sheaves $\cM(-)$ over a subfield $K \subseteq \C$, in the sense of Saito \cite{SaitoFormalism}. We will only use the formal properties familiar from the theory of mixed Hodge modules: the six operations and the weight formalism, together with semisimplicity and decomposition by strict support for pure objects, and for projective direct images, we have the usual purity, the decomposition theorem, and relative Hard Lefschetz results. All objects considered below are of geometric origin. Thus, a reader interested only in mixed Hodge modules may take $K=\C$ and $\mathcal M(X)={\rm MHM}(X,A)$ throughout.

We write $A_X^{\cM}\in D^b\cM(X)$ for the constant object, ${\rm IC}_X^{\cM}$ for the intersection complex, and $(j)$ for the Tate twist. Then, for convenience, we write
\[ \mathbf D_X^{\cM} = \mathbf D(A_X^{\cM}[\dim X])(-\dim X),\]
for the dual constant object. Since the realization functor is conservative, if $X_{\mathbf C}$ is a rational homology manifold of pure dimension $n$, then
\[A_X^{\cM}[n]\cong {\rm IC}_X^{\cM},\]
and in general, we have morphisms 
\[A_X^{\cM}[n] \to {\rm IC}_X^{\cM} \to \mathbf D_X^{\cM}.\]

In any theory of mixed sheaves, we can consider the RHM defect object, the leftmost object in the triangle:
\[ \cK_{X,\cM}^\bullet \to A_{X}^{\cM}[\dim X] \to {\rm IC}_X^{\cM} \xrightarrow[]{+1}.\]

The conclusion of \eqref{eq-RHMDefect} still applies, as Park--Popa's proof does not use any special properties of mixed Hodge modules.

If $\kappa:X\to {\rm Spec}(K)$ is the structure morphism, we set
\[ H^i_{\cM}(X)=\cH^i\kappa_*A_X^{\cM},\]
whose underlying complex of $A$-vector spaces agrees with the singular cohomology $H^i(X^{\rm an}_{\C},A)$.

For a projective rational homology manifold, relative Hard Lefschetz gives the usual primitive decomposition, and we denote the resulting primitive pieces by $H^i_{\cM,\mathrm{prim}}(X)$. More generally, for a pure object $M$, a projective morphism $g$ and a relatively ample line bundle $\ell$, we can define the $\ell$-primitive object for any $a\geq 0$: 
\[ (\cH^{-a}g_* M)_{\rm prim} = \ker(\ell^{a+1} \colon \cH^{-a} g_* M \to \cH^{a+2} g_* M(a+1)). \]

Besides mixed Hodge modules, one may use Saito's theory of ``systems of realizations'' \cite{SaitoCycle} as $\mathcal M(-)$. In this case, the statements below simultaneously retain the corresponding Hodge-theoretic and $\ell$-adic/Galois realizations.

\subsection{One-step resolution}
We now formulate the geometric setting in which we will apply the above formalism. This will be called the ``one-step resolution'' set-up below.

Consider the following Cartesian diagram
\[ \begin{tikzcd} \widetilde{Z} \ar[d,"p"] \ar[r] & \widetilde{X}\ar[d,"f"] \\ Z \ar[r,"i"] & X \end{tikzcd}\]
where $f$ (hence $p$) is projective, $i$ is a closed embedding, and all the varieties are equidimensional. We use the following notation for (co)-dimensions:
\begin{itemize}
    \item $n = \dim X = \dim \widetilde{X}$,
    \item $d_Z = \dim Z$,
    \item $d_{\widetilde{Z}} = \dim \widetilde{Z}$,
    \item $c_Z = n -d_Z$,
    \item $c_{\widetilde{Z}} = n-d_{\widetilde{Z}}$,
    \item $\delta = c_{\widetilde{Z}}-1$,
    \item $d = d_{\widetilde{Z}} - d_Z = c_Z - c_{\widetilde{Z}}$, which we always assume is positive.
\end{itemize}

Moreover, let $L$ be an $f$-relatively ample line bundle on $\widetilde{X}$. We assume that $f_{\C}$ is an isomorphism over the complement of $Z\times_K \C$. 

In this section, we prove the following technical results:
\begin{thm} \label{thm-general} Assume $\widetilde{X} \times_K \C,\widetilde{Z}\times_K \C$ are rational homology manifolds. Then we have isomorphisms for all $0 < j < c_Z-1$:
\[ \cH^{-j}A^{\cM}_X[n] \cong \begin{cases} \bigoplus_{r=0}^{\delta} i_* (\cH^{\delta-j-2r}p_* A_{\widetilde{Z}}^\cM[d_{\widetilde{Z}}])_{\rm prim}(-r) & \delta \leq j \\ \bigoplus_{r = 0}^{j} i_*(\cH^{j-\delta-2r}p_* A_{\widetilde{Z}}^\cM[d_{\widetilde{Z}}])_{\rm prim}(j-\delta-r) & \delta > j\end{cases}\]
and an isomorphism
\[ {\rm Gr}^W_{n-1} \cH^0A^{\cM}_X[n] \cong i_*(\cH^{-\delta}p_* A_{\widetilde{Z}}^\cM[d_{\widetilde{Z}}])_{\rm prim}(-\delta).\]
\end{thm}

With mixed Hodge module coefficients, we can assume partial rational homology manifold conditions and still obtain concrete statements. As the statement concerns a Hodge filtration piece, we assume $X$, hence $Z$, is embedded (locally, if necessary) into some smooth variety $Y$. We will use $i\colon Z \to Y$ also for the closed embedding into $Y$.

\begin{thm}\label{thm-generalMHM}
Let $K = \C$ and assume that ${\rm HRH}(\widetilde{X}) \geq k \geq 0$ and ${\rm HRH}(\widetilde{Z}) \geq k$. Then for $0 < j < c_Z -1$, we have (for underlying filtered \emph{right} $\cD$-modules)
\[F_{k-n} \cH^j \mathbf D_X^H \cong \begin{cases} \bigoplus_{r=0}^{j} F_{k+r-d_{\widetilde{Z}}} i_*(\cH^{j-\delta-2r}p_* {\rm IC}_{\widetilde{Z}}^H)_{\rm prim} & \delta  > j \\ \bigoplus_{r = 0}^{\delta} F_{k+r+j+1-n} i_*(\cH^{\delta-j-2r}p_*{\rm IC}_{\widetilde{Z}}^H)_{\rm prim} & \delta \leq j\end{cases}\]
and 
\[ F_{k-n} {\rm Gr}^W_{n+1} \cH^0\mathbf D_X^H \cong F_{k-d_{\widetilde{Z}}} i_*(\cH^{-\delta}p_* {\rm IC}_{\widetilde{Z}}^H)_{\rm prim} .\]
\end{thm}

The proofs of both results proceed by first expressing the relevant cohomology objects on $X$ as kernels and cokernels of natural restriction maps associated to $f$ and $p$. We then use relative Hard Lefschetz to identify these kernels and cokernels in terms of the primitive pieces of the cohomology of $p$.

As $f_{\C}$ is an isomorphism away from $Z$, we have the dual exact triangles
\begin{equation} \label{eq-tri1} A^{\cM}_X[n] \to f_* A^{\cM}_{\widetilde{X}}[n] \to i_*S \xrightarrow[]{+1}\end{equation}
\begin{equation} \label{eq-tri1DUal} i_*S' \to f_* \mathbf D_{\widetilde{X}}^\cM \to \mathbf D_{X}^\cM \xrightarrow[]{+1}\end{equation}
for some $S,S'\in D^b\cM(Z)$.

Moreover, by adjunction, the morphism
\[ f_* A^{\cM}_{\widetilde{X}}[n] \to i_* S\]
factors as
\[ f_* A^{\cM}_{\widetilde{X}}[n] \to i_*i^* f_* A^{\cM}_{\widetilde{X}}[n] \to  i_* S,\]
and similarly the morphism
\[ i_* S' \to f_* \mathbf D_{\widetilde{X}}^\cM\text{ factors through } i_* i^! f_* \mathbf D_{\widetilde{X}}^\cM.\]

By Saito's base change \cite{SaitoFormalism}*{Prop. 3.10}, we will identify $i^* f_* A^{\cM}_{\widetilde{X}}[n] = p_* A_{\widetilde{Z}}^\cM[n]$ and $i^! f_* \mathbf D_{\widetilde{X}}^{\cM} = p_* \mathbf D_{\widetilde{Z}}^{\cM}(-c_{\widetilde{Z}})[-c_{\widetilde{Z}}]$. Thus, taking cohomology, we have
\[\cH^j i_* i^* f_* A^{\cM}_{\widetilde{X}}[n] = i_* \cH^{c_{\widetilde{Z}}+j} p_*A_{\widetilde{Z}}^\cM[d_{\widetilde{Z}}]\]
\[\cH^j i_* i^! f_* \mathbf D_{\widetilde{X}}^\cM  = i_* \cH^{j-c_{\widetilde{Z}}} p_*\mathbf D_{\widetilde{Z}}^\cM(-c_{\widetilde{Z}})\]

Note that the morphism $f_* A^{\cM}_{\widetilde{X}}[n] \to i_* i^* f_* A^{\cM}_{\widetilde{X}}[n]$ is identified under this base change with the result of applying the functor $f_*$ to the restriction morphism $A^{\cM}_{\widetilde{X}}[n] \to i_{\widetilde{Z},*} A_{\widetilde{Z}}^\cM[n]$, where $i_{\widetilde{Z}} \colon \widetilde{Z} \to \widetilde{X}$ is the inclusion, and the analogous claim holds for the dual.

\begin{lem} \label{lem-computeS} For all $j > -c_Z$, the natural maps
\[ \cH^{c_{\widetilde{Z}}+j} p_* A_{\widetilde{Z}}^\cM[d_{\widetilde{Z}}] \to \cH^j i^* f_* A^{\cM}_{\widetilde{X}}[n] \to \cH^j S\]
are isomorphisms. 

Similarly, for all $j < c_Z$, the maps
\[ \cH^j S' \to \cH^j i^! f_* \mathbf D_{\widetilde{X}}^\cM \to \cH^{j-c_{\widetilde{Z}}} p_* \mathbf D_{\widetilde{Z}}^{\cM}(-c_{\widetilde{Z}})\]
are isomorphisms.
\end{lem}
\begin{proof} We prove the first claim; the second follows by duality (or an analogous argument). The first map is an isomorphism as described above by base change, with no assumption on $j$.

By applying $i^*(-)$ to the triangle \eqref{eq-tri1}, we get the exact triangle
\[A_Z^\cM[d_Z][c_Z] \to i^* f_* A^{\cM}_{\widetilde{X}}[n] \to S \xrightarrow[]{+1}.\]

Now, $\cH^j(A_Z^\cM[d_Z][c_Z]) = 0$ for all $j > -c_Z$, or even $j\neq -c_Z$ if $Z$ is CCI, which proves the claim by looking at the long exact sequence in cohomology.
\end{proof}

\begin{cor} \label{cor-comparefstar} For all $j > 0$, the natural maps
\[ \cH^j f_* A^{\cM}_{\widetilde{X}}[n] \to \cH^j i_* i^* f_* A^{\cM}_{\widetilde{X}}[n] \to i_* \cH^{c_{\widetilde{Z}}+j} p_* A_{\widetilde{Z}}^\cM[d_{\widetilde{Z}}]\]
are isomorphisms, and for $j < 0$, the maps
\[  i_* \cH^{j-c_{\widetilde{Z}}} p_* \mathbf D_{\widetilde{Z}}^\cM(-c_{\widetilde{Z}}) \to \cH^j i_* i^! f_* \mathbf D_{\widetilde{X}}^{\cM} \to\cH^j f_* \mathbf D^{\cM}_{\widetilde{X}} \]
are isomorphisms.
\end{cor}
\begin{proof} As $A^{\cM}_X[n]$ satisfies $\cH^j A^{\cM}_X[n] = 0$ for all $j>0$, we have isomorphisms
\[ \cH^j f_* A^{\cM}_{\widetilde{X}}[n] \to \cH^j i_* S,\]
and so the claim follows from the previous lemma. The second claim follows by a similar argument (or duality).
\end{proof}

For $j \geq 0$, define \[\alpha_j \colon \cH^{-j}f_* A^{\cM}_{\widetilde{X}}[n] \to i_* \cH^{c_{\widetilde{Z}}-j} p_* A_{\widetilde{Z}}^\cM[d_{\widetilde{Z}}]\] and \[\beta_j \colon i_* \cH^{j-c_{\widetilde{Z}}} p_* \mathbf D_{\widetilde{Z}}^{\cM}(-c_{\widetilde{Z}}) \to \cH^j f_* \mathbf D_{\widetilde{X}}^{\cM}.\]

So using the isomorphisms above, the long exact sequence in cohomology coming from the Triangle \eqref{eq-tri1} can be rewritten as
\begin{equation} \label{eq-LESConstantHM} \dots \xrightarrow[]{\alpha_{j+1}} i_* \cH^{c_{\widetilde{Z}}-j-1} p_*A_{\widetilde{Z}}^\cM[d_{\widetilde{Z}}] \to \cH^{-j}A^{\cM}_X[n] \to \cH^{-j}f_* A^{\cM}_{\widetilde{X}}[n] \xrightarrow[]{\alpha_j} \dots \end{equation}
\begin{equation} \label{eq-LESDualHM} \dots \xrightarrow[]{\beta_{j}}  \cH^{j}f_* \mathbf D^{\cM}_{\widetilde{X}} \to \cH^{j}\mathbf D^{\cM}_X \to i_*\cH^{j+1-c_{\widetilde{Z}}} p_*\mathbf D_{\widetilde{Z}}^\cM(-c_{\widetilde{Z}}) \xrightarrow[]{\beta_{j+1}} \dots \end{equation}
as long as $0 \leq j <c_Z -1$. Indeed, for larger $j$, we cannot necessarily apply \lemmaref{lem-computeS} above.

The technical crux of the argument is to identify the kernels and cokernels of the morphisms $\alpha_j$ and $\beta_j$.

In the mixed sheaves setting, we assume $\widetilde{X} \times_K \C$ is a rational homology manifold, or equivalently, that $A^\cM_{\widetilde{X}}[n]$ (hence, $\mathbf D_{\widetilde{X}}^{\cM}$) is pure. 

\begin{lem} \label{lem-RHMSES} Assume $\widetilde{X}\times_K \C$ is a rational homology manifold. We get short exact sequences for all $0 < j < c_Z-1$:
\[ 0 \to \cH^{-j-1} f_* A^{\cM}_{\widetilde{X}}[n] \xrightarrow[]{\alpha_{j+1}} i_*\cH^{c_{\widetilde{Z}}-j-1} p_* A_{\widetilde{Z}}^\cM[d_{\widetilde{Z}}] \to \cH^{-j}A^{\cM}_X[n] \to 0\]
and short exact sequences
\[ 0 \to \cH^{-1} f_* A^{\cM}_{\widetilde{X}}[n] \xrightarrow[]{\alpha_{1}} i_*\cH^{c_{\widetilde{Z}}-1} p_* A_{\widetilde{Z}}^\cM[d_{\widetilde{Z}}] \to W_{n-1}\cH^{0}A^{\cM}_X[n] \to 0,\]
\[ 0 \to {\rm IC}_X^{\cM}  \to \cH^{0} f_* A^{\cM}_{\widetilde{X}}[n] \xrightarrow[]{\alpha_{0}} i_*\cH^{c_{\widetilde{Z}}} p_* A_{\widetilde{Z}}^\cM[d_{\widetilde{Z}}] \to 0.\]

Dually, we have short exact sequences for all $0 < j < c_Z-1$:
\[ 0 \to \cH^{j}\mathbf D^{\cM}_X \to i_*\cH^{j+1-c_{\widetilde{Z}}} p_*\mathbf D_{\widetilde{Z}}^\cM(-c_{\widetilde{Z}}) \xrightarrow[]{\beta_{j+1}} \cH^{j+1} f_* \mathbf D_{\widetilde{X}}^{\cM} \to 0\]
and short exact sequences
\[ 0 \to \cH^0\mathbf D_X^\cM/{\rm IC}_X^{\cM} \to i_*\cH^{1-c_{\widetilde{Z}}} p_*\mathbf D_{\widetilde{Z}}^\cM(-c_{\widetilde{Z}}) \xrightarrow[]{\beta_{1}} \cH^{1} f_* \mathbf D_{\widetilde{X}}^{\cM} \to 0\]
\[ 0 \to i_*\cH^{-c_{\widetilde{Z}}} p_*\mathbf D_{\widetilde{Z}}^\cM(-c_{\widetilde{Z}}) \xrightarrow[]{\beta_{0 }} \cH^{0} f_* \mathbf D_{\widetilde{X}}^{\cM} \to {\rm IC}_X^{\cM} \to 0.\]

\end{lem}
\begin{proof} We know that for $j>0$, the object $\cH^{-j}A^{\cM}_X[n]$ has weights $\leq n-j-1$ by \eqref{eq-RHMDefect}. Thus, because $\cH^{-j} f_* A^{\cM}_{\widetilde{X}}[n]$ has weight $n-j$, the morphism
\[ \cH^{-j}A^{\cM}_X[n] \to \cH^{-j} f_* A^{\cM}_{\widetilde{X}}[n]\]
is zero, proving the first claim.

The other two exact sequences follow similarly, by taking $W_{n-1}(-)$ or ${\rm Gr}^W_{n}(-)$ of the long exact sequence above.

The same argument (or duality) gives the claims about the dual objects.
\end{proof}

In the mixed Hodge module setting, we can give a similar result with a weaker assumption. Working locally, we can assume $X$ (hence $Z$) can be embedded into a smooth variety, and as $f$ is projective, the same is true for $\widetilde{X}$ (hence $\widetilde{Z}$). Then we can discuss the underlying filtered $\cD$-module for any mixed Hodge module. In the statements below, we use right filtered $\cD$-modules unless otherwise specified.

We will make use of the following vanishing theorem for certain Hodge filtered pieces of direct images.

\begin{lem} \label{lem-VanishingPushforwardHodge} Let $M^\bullet \in D^{b}({\rm MHM}(Y))$ and let $\pi \colon Y \to Y'$ be a projective morphism between smooth varieties. Let $(\cM^\bullet,F)$ be the underlying complex of filtered right $\cD_Y$-modules and assume $F_p \cM^\bullet = 0$. Then we have
\[ F_{p} \cH^j \pi_*(\cM^\bullet) = 0\text{ for all } j \in \Z.\]

If $M^\bullet \in D^{\geq 0}({\rm MHM}(Y))$, then we get additional vanishing for the negative cohomology modules:
\[ F_{p+j} \cH^{-j} \pi_*(\cM^\bullet) = 0 \text{ for all } j \geq 0.\]
\end{lem}
\begin{proof} Let $(\cM^\bullet,F)$ denote the underlying complex of right $\cD_Y$-modules. We first study behavior for the individual cohomology modules $\cH^j M^\bullet$, so let $M$ be a single mixed Hodge module on $Y$ with $F_p \cM = 0$. We will prove the desired vanishings for the module $M$. 

We have the graph embedding $\Gamma \colon Y \to Y' \times \P^N$ and the projection $\Pi \colon Y' \times \P^N \to Y'$, so that $\pi = \Pi \circ \Gamma$ and thus it suffices to replace $M$ by $\Gamma_* M$ and $\pi$ by $\Pi$. By strictness of pushforwards, we have for all $q\in \Z$ the equality
\[ F_q \cH^j \pi_*(\cM) = R^j \pi_*( F_q {\rm DR}_{Y' \times \P^N/Y'}(\cM)),\]
where ${\rm DR}_{Y'\times \P^N/Y'}(-)$ is the relative de Rham complex.

The filtration on the relative de Rham complex is
\[ F_{p-N} \cM \otimes \bigwedge^N \cT_{\P^N} \to F_{p-(N-1)} \cM \otimes \bigwedge^{N-1} \cT_{\P^N} \to \dots \to F_p \cM,\]
placed in degrees $-N,\dots, 0$, and so we see that
\[ \cH^b F_q {\rm DR}_{Y'\times \P^N/Y'}(\cM) = 0 \text{ if } q+b \leq p,\]
and so we see that $F_p {\rm DR}_{Y'\times \P^N/Y'}(\cM) = 0$, proving the first claim.

For the second claim, we have the spectral sequence
\[ E_2^{a,b} = R^a \pi_*( \cH^b F_q {\rm DR}_{Y' \times \P^N/Y'}(\cM)) \implies R^{a+b} \pi_*(F_q {\rm DR}_{Y' \times \P^N/Y'}(\cM)),\]
and so for $q$ fixed, we see that $E_2^{a,b} = 0$ for $a < 0$, for $b \notin [-N,0]$ or for $q + b  \leq p$. For $j > 0$ fixed, we consider all $a+b = -j$ with $a \geq 0$ and $b \in [-N,0]$. In particular, we really have $b \in [-N,-j]$. But then $E_2^{a,b} = 0$ for all $q+b \leq p$. So all relevant $E_2$-page pieces vanish if we take $q = p+j$, so that $q+b = (p+j) +b \leq p$ for all $b \in [-N,-j]$. This proves the second claim.

Now, we consider an arbitrary $M^\bullet \in D^{b}({\rm MHM}(Y))$. Consider the spectral sequence $E_2^{i,j} = \cH^i \pi_*(\cH^j(M^\bullet)) \implies \cH^{i+j} \pi_*(M^\bullet)$. By the above case, we see that $F_p E_2^{i,j} = 0$ for all $i,j$ and
\[ F_{p-i} E_2^{i,j} =0\text{ for all } i \leq 0.\]

Thus, the same is true for the $E_\infty$ pages (as the differentials are strict with respect to the Hodge filtration). So we immediately see that the first claim holds.

For the second claim, assume $M^\bullet \in D^{\geq 0}({\rm MHM}(Y))$ and let $a > 0$ and assume $i+j = -a$. Then since $j \geq 0$ we see that $i \leq -a$. Thus, for all relevant $E_\infty$ terms, we see
\[ F_{p+a} E_\infty^{i,j} \subseteq F_{p-i} E_{\infty}^{i,j} = 0,\]
which proves the claim.
\end{proof}

\begin{lem} \label{lem-RHMSESMHM} Assume ${\rm HRH}(\widetilde{X}) \geq k$ for some $k\in \Z_{\geq 0}$. We get short exact sequences for all $0 < j < c_Z-1$:
\[ 0 \to F_{k-n}\cH^{-j-1} f_* \Q^H_{\widetilde{X}}[n] \xrightarrow[]{\alpha_{j+1}} F_{k-n}i_*\cH^{c_{\widetilde{Z}}-j-1} p_* \Q^H_{\widetilde{Z}}[d_{\widetilde{Z}}] \to F_{k-n} \cH^{-j}\Q^H_X[n] \to 0\]
and the short exact sequence
\[ 0 \to F_{k-n}\cH^{-1} f_* \Q^H_{\widetilde{X}}[n] \xrightarrow[]{\alpha_{1}} F_{k-n}i_*\cH^{c_{\widetilde{Z}}-1} p_* \Q^H_{\widetilde{Z}}[d_{\widetilde{Z}}] \to F_{k-n} W_{n-1}\cH^{0}\Q^H_X[n] \to 0.\]

Dually, we have short exact sequences for all $0 < j < c_Z-1$:
\[ 0 \to F_{k-n}\cH^{j}\mathbf D^{H}_X \to F_{k-d_{\widetilde{Z}}}i_*\cH^{j+1-c_{\widetilde{Z}}} p_*\mathbf D_{\widetilde{Z}}^H \xrightarrow[]{\beta_{j+1}} F_{k-n}\cH^{j+1} f_* \mathbf D_{\widetilde{X}}^{H} \to 0\]
and the short exact sequence
\[ 0 \to F_{k-n}(\cH^0\mathbf D_X^H/{\rm IC}_X^H) \to F_{k-d_{\widetilde{Z}}} i_*\cH^{1-c_{\widetilde{Z}}} p_*\mathbf D_{\widetilde{Z}}^H \xrightarrow[]{\beta_{1}} F_{k-n}\cH^{1} f_* \mathbf D_{\widetilde{X}}^{H} \to 0.\]
\end{lem}
\begin{proof} Recall \cite{DOR1}*{Thm. D(1)} that ${\rm HRH}(\widetilde{X}) \geq k$ if and only if $F_{k-n} \gamma^\vee \colon F_{k-n} {\rm IC}_{\widetilde{X}}^{H} \to F_{k-n} \mathbf D_{\widetilde{X}}^H$ is a quasi-isomorphism and this implies that the map $F_{k-n} \gamma \colon F_{k-n} \Q_{\widetilde{X}}^H[n] \to F_{k-n} {\rm IC}_{\widetilde{X}}^H$ is also a quasi-isomorphism.

Thus, by the first claim of \lemmaref{lem-VanishingPushforwardHodge} applied to the cones of the morphisms $\gamma$ and $\gamma^\vee$ (using the fact that $\dim \widetilde{X} = \dim X$), we see that the natural morphisms
\[  F_{k-n} \cH^j f_* \Q_{\widetilde{X}}^H[n] \to F_{k-n} \cH^j f_* {\rm IC}_{\widetilde{X}}^H\]
\[ F_{k-n} \cH^j f_* {\rm IC}_{\widetilde{X}}^H \to F_{k-n} \cH^j f_* \mathbf D_{\widetilde{X}}^H\]
are isomorphisms for all $j \in \Z$. As $f$ is projective and ${\rm IC}_{\widetilde{X}}^H$ is pure of weight $n$, we know that $\cH^j f_* {\rm IC}_{\widetilde{X}}^H$ is pure of weight $n+j$ for all $j\in \Z$. Again, using \eqref{eq-RHMDefect}, we see that the morphisms $\cH^{-j} \Q_X^H[n] \to \cH^{-j} f_* {\rm IC}_{\widetilde{X}}^H$ and $\cH^j f_* {\rm IC}_{\widetilde{X}}^H \to \cH^j \mathbf D_X^H$ are $0$ for $j>0$. We have the commutative triangles
\[ \begin{tikzcd} F_{k-n} \cH^{-j}\Q_X^H[n] \ar[r] \ar[dr,"0"] & F_{k-n} \cH^{-j} f_* \Q_{\widetilde{X}}^H[n] \ar[d,"\cong"] \\ {} & F_{k-n} \cH^{-j} f_* {\rm IC}_{\widetilde{X}}^H \end{tikzcd}, \, \begin{tikzcd} F_{k-n} \cH^{j} f_* {\rm IC}_{\widetilde{X}}^H \ar[d,"\cong"] \ar[dr,"0"] & {} \\  F_{k-n} \cH^{j} f_* \mathbf D_{\widetilde{X}}^H \ar[r] & F_{k-n} \cH^{j}\mathbf D_X^H \end{tikzcd},\]
and so once again the morphisms in the long exact sequence vanish, giving the claimed short exact sequences.
\end{proof}

In the following diagrams and formulas, we suppress the closed-embedding pushforward $i_*$ on objects supported on $Z$.

Back to the general setting of mixed sheaves, for all $j>0$, the following diagram commutes, because the horizontal morphisms are the restriction maps:
\begin{equation} \label{sq-alphaj} \begin{tikzcd} \cH^{-j} f_* A^{\cM}_{\widetilde{X}}[n] \ar[r,"\alpha_j"] \ar[d,"c_1(L)^j"] & \cH^{c_{\widetilde{Z}}-j} p_* A_{\widetilde{Z}}^\cM[d_{\widetilde{Z}}]\ar[d,"c_1(L\vert_{\widetilde{Z}})^j"] \\ \cH^jf_* A^{\cM}_{\widetilde{X}}[n](j) \ar[r,"\cong"]  & \cH^{c_{\widetilde{Z}}+j} p_* A_{\widetilde{Z}}^\cM[d_{\widetilde{Z}}](j) \end{tikzcd},\end{equation}
where the bottom horizontal morphism is the isomorphism from \corollaryref{cor-comparefstar}.  Similarly, the diagram
\begin{equation} \label{sq-betaj} \begin{tikzcd}  \cH^{-j-c_{\widetilde{Z}}} p_* \mathbf D_{\widetilde{Z}}^\cM(-c_{\widetilde{Z}}-j)  \ar[r,"\cong"]   \ar[d,"c_1(L\vert_{\widetilde{Z}})^j"]  & \cH^{-j} f_* \mathbf D^{\cM}_{\widetilde{X}}(-j)\ar[d,"c_1(L)^j"] \\ \cH^{j-c_{\widetilde{Z}}} p_* \mathbf D_{\widetilde{Z}}^\cM(-c_{\widetilde{Z}}) \ar[r,"\beta_j"] &\cH^{j} f_* \mathbf D^{\cM}_{\widetilde{X}} \end{tikzcd},\end{equation}
commutes, with top isomorphism given by \corollaryref{cor-comparefstar}.

If we assume $\widetilde{X}\times_K \C$ is a rational homology manifold, then the left and bottom morphisms in the Diagram \eqref{sq-alphaj} are isomorphisms and therefore induce an isomorphism 
\[ {\rm coker}(\alpha_j) \cong \ker(c_1(L\vert_{\widetilde{Z}})^j),\]
and similarly, as the top and rightmost morphisms in Diagram \eqref{sq-betaj} are isomorphisms, we get an isomorphism
\[ \ker(\beta_j) \cong {\rm coker}(c_1(L\vert_{\widetilde{Z}})^j).\]

Now we assume $\widetilde{Z} \times_K \C$ is a rational homology manifold, or equivalently, we assume $A^{\cM}_{\widetilde{Z}}[d_{\widetilde{Z}}]$ is pure of weight $d_{\widetilde{Z}}$. By \lemmaref{lem-RHMSES}, this immediately implies the following:
\[ \cH^{-j}A^{\cM}_X[n] \text{ is pure of weight } n-j-1 \text{ for all } 0 < j < c_Z -1\]
\[ {\rm Gr}^W_i \cH^0 A^{\cM}_X[n]\neq 0 \implies i \in \{n-1,n\},\]
so that $W_{n-1} \cH^0 A^{\cM}_X[n] = {\rm Gr}^W_{n-1} \cH^0 A^{\cM}_X[n]$ and $\cH^0\mathbf D_X^{\cM}/{\rm IC}_X^{\cM} = {\rm Gr}^W_{n+1} \cH^0\mathbf D_X^{\cM}$.

A simple computation concerning the Lefschetz decomposition in this setting gives the isomorphisms
\begin{equation} \label{eq-identifykerc1j} {\rm coker}(\alpha_j) \cong \ker(c_1(L\vert_{\widetilde{Z}})^j) \cong \begin{cases} \bigoplus_{r=0}^{c_{\widetilde{Z}}-1} (\cH^{c_{\widetilde{Z}}-j-2r}p_* A_{\widetilde{Z}}^\cM[d_{\widetilde{Z}}])_{\rm prim}(-r) & c_{\widetilde{Z}} \leq j \\ \bigoplus_{r = 0}^{j-1} (\cH^{j-c_{\widetilde{Z}}-2r}p_* A_{\widetilde{Z}}^\cM[d_{\widetilde{Z}}])_{\rm prim}(j-c_{\widetilde{Z}}-r) & c_{\widetilde{Z}} >j\end{cases},\end{equation}
\begin{equation} \label{eq-identifycokerc1j} \ker(\beta_j) \cong {\rm coker}(c_1(L\vert_{\widetilde{Z}})^j) \cong \begin{cases} \bigoplus_{r=0}^{j-1} (\cH^{j-c_{\widetilde{Z}}-2r}p_* \mathbf D_{\widetilde{Z}}^\cM)_{\rm prim}(-r) & c_{\widetilde{Z}} > j \\ \bigoplus_{r = 0}^{c_{\widetilde{Z}}-1} (\cH^{c_{\widetilde{Z}}-j-2r}p_* \mathbf D_{\widetilde{Z}}^\cM)_{\rm prim}(c_{\widetilde{Z}}-j-r) & c_{\widetilde{Z}} \leq j\end{cases}.\end{equation}

In the mixed Hodge modules setting, if we assume ${\rm HRH}(\widetilde{X}) \geq k$ (and no longer that $\widetilde{X}$ is a rational homology manifold), then we consider the square for $j>0$:
\begin{equation} \label{sq-betajMHM} \begin{tikzcd}  F_{k+j-d_{\widetilde{Z}}}\cH^{-j-c_{\widetilde{Z}}} p_* \mathbf D_{\widetilde{Z}}^H  \ar[r,"\cong"]   \ar[d,"c_1(L\vert_{\widetilde{Z}})^j"]  & F_{k+j-n} \cH^{-j} f_* \mathbf D^{H}_{\widetilde{X}}\ar[d,"c_1(L)^j"] \\ F_{k-d_{\widetilde{Z}}} \cH^{j-c_{\widetilde{Z}}} p_* \mathbf D_{\widetilde{Z}}^H\ar[r,"\beta_j"] &F_{k-n} \cH^{j} f_* \mathbf D^{H}_{\widetilde{X}} \end{tikzcd},\end{equation}
where the rightmost vertical map is an isomorphism by the assumption on $\widetilde{X}$. Indeed, we have the commutative diagram
\[\begin{tikzcd}   F_{k+j-n}\cH^{-j} f_*{\rm IC}_{\widetilde{X}}^H  \ar[r,"\cong"]   \ar[d,"c_1(L)^j"]  & F_{k+j-n}\cH^{-j} f_* \mathbf D_{\widetilde{X}}^H    \ar[d,"c_1(L)^j"]  \\ F_{k-n} \cH^{j} f_* {\rm IC}_{\widetilde{X}}^H\ar[r,"\cong"] & F_{k-n} \cH^{j} f_* \mathbf D_{\widetilde{X}}^H \end{tikzcd}\]
where the assumption ${\rm HRH}(\widetilde{X}) \geq k$ implies the horizontal morphisms are isomorphisms by \lemmaref{lem-VanishingPushforwardHodge}. The leftmost vertical map is an isomorphism by relative Hard Lefschetz for the pure Hodge module ${\rm IC}_{\widetilde{X}}^H$.

Thus, by the commutative diagram \eqref{sq-betajMHM}, we get an isomorphism
\[ F_{k-n} \ker(\beta_j) \cong F_{k-n}\left({\rm coker}(c_1(L\vert_{\widetilde{Z}})^j)(-c_{\widetilde{Z}})\right).\]

Assume moreover that ${\rm HRH}(\widetilde{Z}) \geq k$. Equivalently, assume the morphism $F_{k-d_{\widetilde{Z}}}{\rm IC}^H_{\widetilde{Z}} \to F_{k-d_{\widetilde{Z}}} \mathbf D^H_{\widetilde{Z}}$ is a quasi-isomorphism. Applying \lemmaref{lem-VanishingPushforwardHodge} to the cone of this morphism, we conclude
\[ F_{k-d_{\widetilde{Z}}} \cH^j p_*{\rm IC}^H_{\widetilde{Z}} \to F_{k-d_{\widetilde{Z}}} \cH^j p_* \mathbf D^H_{\widetilde{Z}}\]
is an isomorphism for all $j\in \Z$, and for $j > 0$, we have that the natural map
\[ F_{k+j-d_{\widetilde{Z}}} \cH^{-j} p_* {\rm IC}^H_{\widetilde{Z}} \to F_{k+j-d_{\widetilde{Z}}} \cH^{-j} p_* \mathbf D^H_{\widetilde{Z}}\]
is an isomorphism. These isomorphisms commute with cupping with $c_1(L\vert_{\widetilde{Z}})^j$.

Thus, if we consider the commutative diagram:
\begin{equation} \label{sq-MHMLefschetz} \begin{tikzcd}   F_{k+j-d_{\widetilde{Z}}}\cH^{-j-c_{\widetilde{Z}}} p_*{\rm IC}_{\widetilde{Z}}^H  \ar[r,"\cong"]   \ar[d,"c_1(L\vert_{\widetilde{Z}})^j"]  & F_{k+j-d_{\widetilde{Z}}}\cH^{-j-c_{\widetilde{Z}}} p_* \mathbf D_{\widetilde{Z}}^H    \ar[d,"c_1(L\vert_{\widetilde{Z}})^j"]  \\ F_{k-d_{\widetilde{Z}}} \cH^{j-c_{\widetilde{Z}}} p_* {\rm IC}_{\widetilde{Z}}^H\ar[r,"\cong"] & F_{k-d_{\widetilde{Z}}} \cH^{j-c_{\widetilde{Z}}} p_* \mathbf D_{\widetilde{Z}}^H \end{tikzcd},\end{equation}
we conclude that the cokernels of the vertical morphisms are isomorphic. The object $p_* {\rm IC}_{\widetilde{Z}}^H$ satisfies the relative Hard Lefschetz theorem, and so, similarly to the computation of \eqref{eq-identifycokerc1j} above, we can write
\[ F_{k-n} \ker(\beta_j) = F_{k-n+c_{\widetilde{Z}}}{\rm coker}(c_1(L\vert_{\widetilde{Z}})^j) = \begin{cases} \bigoplus_{r=0}^{j-1} F_{k+r-d_{\widetilde{Z}}} (\cH^{j-c_{\widetilde{Z}}-2r}p_* {\rm IC}_{\widetilde{Z}}^H)_{\rm prim} & c_{\widetilde{Z}} > j \\ \bigoplus_{r = 0}^{c_{\widetilde{Z}}-1} F_{k+r+j-n} (\cH^{c_{\widetilde{Z}}-j-2r}p_*{\rm IC}_{\widetilde{Z}}^H)_{\rm prim} & c_{\widetilde{Z}} \leq j\end{cases}.\]

Recalling that $\delta = c_{\widetilde{Z}}-1$ and applying the short exact sequences in \lemmaref{lem-RHMSES} (and \lemmaref{lem-RHMSESMHM} in the mixed Hodge module case) this proves \theoremref{thm-general} and \theoremref{thm-generalMHM}.

At this point, there are two natural situations one could consider. In future work, the authors will focus on the following situation: the variety $\widetilde{X} \times_K \C$ is a rational homology manifold, the morphism $p \colon \widetilde{Z} \to Z$ is a smooth map between smooth varieties of positive dimension.

We focus on the case where $Z$ is a point and $\widetilde{X} \times_K \C,\widetilde{Z}\times_K \C$ are rational homology manifolds. We let $d =d_{\widetilde{Z}}$ to simplify the notation below (this agrees with our earlier definition of $d$, as we now assume $d_Z = 0$). 

\begin{cor} \label{cor-Isolated} Assume $Z$ is a point and that $\widetilde{X} \times_K \C,\widetilde{Z} \times_K \C$ are irreducible rational homology manifolds. Assume $\dim X \geq 2$.

Then
\begin{enumerate} \item ${\rm lcdef}(X) \leq \dim X -2$.
    \item For all $0 < j \leq \dim X -2$, we have an isomorphism of pure mixed sheaves of weight $\dim X+j+1$:
    \[ \cH^{-j} \cK_{\cM,X}^\bullet(-j-1) \cong \cH^j \mathbf D^\cM_X \cong \begin{cases} \bigoplus_{r=0}^{\delta} i_*H_{\cM,\rm prim}^{d-(j-\delta+2r)}(\widetilde{Z})(-j-r-1) & \delta \leq j \\ \bigoplus_{r = 0}^{j} i_*H_{\cM,\rm prim}^{d-(\delta-j+2r)}(\widetilde{Z})(-\delta-r-1) & \delta >j \end{cases}.\]
    \item We have ${\rm Gr}^W_i \cH^0\mathbf D_X^\cM \neq 0$ implies $i \in \{\dim X,\dim X+1\}$, and an isomorphism
    \[ \cH^0 \cK_{\cM,X}^\bullet (-1) \cong {\rm Gr}^W_{\dim X +1} \cH^0 \mathbf D_X^\cM \cong i_* H^{d-\delta}_{\cM,\rm prim}(\widetilde{Z})(-\delta-1).\]

\item We have an isomorphism in $D^b(\cM(X))$ 
\[ f_* A^{\cM}_{\widetilde{X}}[\dim X] \cong {\rm IC}_X^\cM \oplus  i_* H^{\dim X}_{\cM}(\widetilde{Z}) \oplus \bigoplus_{\ell =1}^{d-c_{\widetilde{Z}}} \left(i_*H^{\dim X +\ell}_{\cM}(\widetilde{Z})[-\ell] \oplus i_*H^{\dim X+\ell}_{\cM}(\widetilde{Z})(\ell)[\ell]\right).\]
\end{enumerate}
\end{cor}
\begin{proof} We prove the claim about the local cohomological defect, which is really only a claim about the underlying $A$-perverse sheaves. We have the exact triangle
\[A^{\cM}_X[n] \to f_* A^{\cM}_{\widetilde{X}}[n] \to i_*S \xrightarrow[]{+1}\]
and by applying $i^*$, we get
\[ A^\cM[n] \to p_* A_{\widetilde{Z}}^\cM[n] \to S \xrightarrow[]{+1}.\]

By the long exact sequence in cohomology, we get
\[ 0 \to \cH^{-n-1} p_* A_{\widetilde{Z}}^\cM[n] \to \cH^{-n-1} S \to A^\cM \xrightarrow[]{\chi} \cH^{-n}p_* A_{\widetilde{Z}}^\cM[n] \to \cH^{-n} S \to 0,\]
but we have $\cH^{-n-1}p_* A_{\widetilde{Z}}^\cM[n] = 0$ and $\cH^{-n}p_* A_{\widetilde{Z}}^\cM[n] = H^0_{\cM}(\widetilde{Z})$ by comparing with the underlying $A$-perverse sheaves, so that the map $\chi$ is an isomorphism. This proves $\cH^{-n-1} S = \cH^{-n} S = 0$. Similarly, we have isomorphisms $0 = \cH^{-n-j} p_* A_{\widetilde{Z}}^\cM[n] \to \cH^{-n-j} S$ for all $j \geq 2$. This implies the same vanishing holds for $i_* S$.

We get the exact sequence
\[ 0 \to \cH^{-n+1} A^{\cM}_X[n] \to \cH^{-n+1} f_* A^{\cM}_{\widetilde{X}}[n] \xrightarrow[]{\alpha_{n-1}} \cH^{-n+1} i_* S,\]
but we argued above that $\alpha_{n-1}$ is injective. Thus, $\cH^{-n+1}A^{\cM}_X[n] = 0$.

Using that $f$ is an isomorphism away from the point $Z$, its defect of semismallness is easily computed to be $\max\{0, 2d-n\}$. Thus, $\cH^{-j} f_* A^{\cM}_{\widetilde{X}}[n] = 0$ for all $j > \max\{0, 2d-n\}$. We can write $2d-n = n - 2c_{\widetilde{Z}}$. As $c_{\widetilde{Z}} \geq 1$, we see then that for all $j\geq \dim X$, we have
\[ \cH^{-j} A^{\cM}_X[n] \cong \cH^{-j}f_* A^{\cM}_{\widetilde{X}}[n] = 0,\]
which proves the bound on the local cohomological defect.

The computation of $\cH^j \mathbf D_X^\cM$ and ${\rm Gr}^W_{n+1} \cH^0\mathbf D_X^\cM$ follows by \theoremref{thm-general} and duality.

For the last claim, we apply the Decomposition theorem to the pure complex $f_* A^{\cM}_{\widetilde{X}}[n]$, giving an isomorphism
\[ f_* A^{\cM}_{\widetilde{X}}[n] \cong \bigoplus_{\ell \in \Z} \cH^\ell(f_* A^{\cM}_{\widetilde{X}}[n])[-\ell].\]

We have the short exact sequence
\[ 0 \to {\rm IC}_X^\cM \to \cH^0f_* A^{\cM}_{\widetilde{X}}[n] \xrightarrow[]{\alpha_0} i_*H^n_{\cM}(\widetilde{Z}) \to 0,\]
and for $\ell > 0$, we have the isomorphism
\[ \cH^\ell f_* A^{\cM}_{\widetilde{X}}[n] \to i_* \cH^{c_{\widetilde{Z}}+\ell} p_* A^{\cM}_{\widetilde{Z}}[d] = i_* H^{n+\ell}_{\cM}(\widetilde{Z}),\]
by \corollaryref{cor-comparefstar}. For the negative cohomologies, we use the Lefschetz isomorphisms for $f_* A^{\cM}_{\widetilde{X}}[n]$ to conclude that, for $\ell >0$, we have an isomorphism
\[\cH^{-\ell}f_* A^{\cM}_{\widetilde{X}}[n] \cong \cH^\ell f_* A^{\cM}_{\widetilde{X}}[n] (\ell) = i_* H^{n+\ell}_{\cM}(\widetilde{Z})(\ell),\]
as claimed.
\end{proof}

\begin{rmk} \label{rmk-IHIsolated} Assume $\widetilde{Z}\times_K \C$ and $\widetilde{X} \times_K \C$ are irreducible rational homology manifolds. We have the following isomorphisms for all $j\in \Z$:
\[ H^j_{\cM}(\widetilde{X}) \cong {\rm IH}^j_{\cM}(X) \oplus \begin{cases} H^{j}_{\cM}(\widetilde{Z}) & j \geq n \\ H^{2n-j}_{\cM}(\widetilde{Z})(n-j) & j < n \end{cases}.\]   
\[ H^j_{\cM,c}(\widetilde{X}) \cong {\rm IH}^j_{\cM,c}(X) \oplus \begin{cases} H^{j}_{\cM}(\widetilde{Z}) & j \geq n \\ H^{2n-j}_{\cM}(\widetilde{Z})(n-j) & j < n \end{cases}.\]
\end{rmk}

Under the weaker assumption (in the mixed Hodge modules setting) we can understand the Hodge filtration on the cohomology of $\mathbf D_X^H$, which is the main object of interest. Let $X \subseteq Y$ be an embedding into a smooth variety and let $Z$ be defined by the vanishing of some coordinates $y_1,\dots, y_N$ on $Y$. By choosing such coordinates, the pushforward along $i \colon Z \to Y$ has a simple form. Then
\begin{cor} \label{cor-HodgeFiltDualObject} Assume $Z$ is a point (defined by $y_1,\dots, y_N$ in the smooth variety $Y$), ${\rm HRH}(\widetilde{X}) \geq k \geq 0$ and ${\rm HRH}(\widetilde{Z}) \geq k$. Then for all $0 < j \leq \dim X-2$ we have
\[ F_{k-n} \cH^j \mathbf D_X^H = \begin{cases} \bigoplus_{r=0}^j \left( \bigoplus_{\alpha \in \N^N} F_{k+r-d-|\alpha|} {\rm IH}^{d-(\delta -j +2r)}_{\rm prim}(\widetilde{Z})\de_y^\alpha\right) & \delta > j \\  \bigoplus_{r=0}^{\delta} \left(\bigoplus_{\alpha \in \N^N} F_{k+r+j+1-n-|\alpha|} {\rm IH}^{d-(j-\delta +2r)}_{\rm prim}(\widetilde{Z})\de_y^\alpha\right) & \delta \leq j \end{cases}\]
and
\[ F_{k-n}( \cH^0\mathbf D_X^H/{\rm IC}_X^H) = \left(\bigoplus_{\alpha \in \N^N} F_{k-d-|\alpha|} {\rm IH}^{d-\delta}_{\rm prim}(\widetilde{Z})\de_y^\alpha\right).\]

In particular, we have
\[ F_{k-n} \cH^j \mathbf D_X^H = 0 \text{ if and only if } \begin{cases} \bigoplus_{r=0}^j F_{k+r-d} {\rm IH}^{d-(\delta -j +2r)}_{\rm prim}(\widetilde{Z}) = 0 & \delta > j \\  \bigoplus_{r=0}^{\delta} F_{k+r+j+1-n} {\rm IH}^{d-(j-\delta +2r)}_{\rm prim}(\widetilde{Z}) = 0 & \delta \leq j \end{cases}\]
and
\[ F_{k-n}( \cH^0\mathbf D_X^H/{\rm IC}_X^H) = 0 \text{ if and only if } F_{k-d} {\rm IH}^{d-\delta}_{\rm prim}(\widetilde{Z}) = 0.\]
\end{cor}

\subsection{Cohomology Computations}
Under the assumptions of \corollaryref{cor-Isolated}, assume further that $\delta = 0$, meaning $\widetilde{Z}$ has codimension $1$ in $\widetilde{X}$. Then we give a concrete description of the cohomology $H^\bullet_{\cM}(X)$ in any mixed sheaf category, under the assumption that $X$ is proper.

Indeed, recall the defining triangle for the RHM defect object:
\[ \cK_{X,\cM}^\bullet \to A_X^{\cM}[\dim X] \to {\rm IC}_X^{\cM} \xrightarrow[]{+1}.\] Note that \theoremref{thm-general} proves that $\cK_{X,\cM}^\bullet$ is a pure object (in the derived sense) of weight $\dim X -1$. If we shift,
\begin{equation} \label{eq-triangleShifted} \cK_{X,\cM}^\bullet[-\dim X] \to A_X^{\cM} \to {\rm IC}_X^{\cM}[-\dim X] \xrightarrow[]{+1},\end{equation}
then the leftmost term is pure of weight $-1$ and the rightmost term is pure of weight $0$. In other words, for all $j\in \Z$, 
\[ \cH^j\kappa_*(\cK_{X,\cM}^\bullet[-\dim X]) \text{ is pure of weight } j-1,\]
and as is well known, 
\[ \cH^j \kappa_* ({\rm IC}_X^{\cM}[-\dim X]) = {\rm IH}^j_{\cM}(X) \text{ is pure of weight } j.\]

As $\cK_{X,\cM}^\bullet$ is supported on a point, we can easily compute from \theoremref{thm-general}
\[ \cH^{-j} \kappa_*(\cK^\bullet_{X,\cM}) \cong H^{d-j}_{\cM, {\rm prim}}(\widetilde{Z}).\]

\begin{lem} In our setting, for all $j > 0$, we have
\[ \cH^{j+1}\kappa_*(\cK_{X,\cM}^\bullet[-\dim X]) \cong {\rm IH}^{j}_{x,\cM}(X) = \cH^{j}(i_x^* {\rm IC}_X^{\cM}[-\dim X]),\]
where $i_x\colon \{x\} \to X$ is the inclusion, and the map
\[ {\rm IH}^j_{\cM}(X) \to \cH^{j+1} \kappa_*(\cK_{X,\cM}^\bullet[-\dim X])\]
is identified with the natural map
\[ {\rm IH}^j_{\cM}(X) \to {\rm IH}^j_{x,\cM}(X).\]
\end{lem}
\begin{proof} By assumption, $\cK_{X,\cM}^\bullet$ is supported on $\{x\} \subseteq X$. Thus, the natural map
\[ {\rm IC}_X^{\cM}[-\dim X] \to \cK_{X,\cM}^\bullet[-\dim X+1]\]
factors, by adjunction, through
\[ i_{x,*} i_x^* {\rm IC}_X^{\cM}[-\dim X].\]

Then, by applying $i_x^*$ to a shift of \eqref{eq-triangleShifted}, we get the triangle
\[ A^{\cM} \to i_x^* {\rm IC}_X^{\cM}[-\dim X] \to i_x^* \cK_{X,\cM}^\bullet[-\dim X+1] \xrightarrow[]{+1},\]
and so for $j > 0$, we have an isomorphism
\[ \cH^j i_x^*{\rm IC}_X^{\cM}[-\dim X] \cong \cH^{j} i_x^* \cK_{X,\cM}^\bullet[-\dim X+1].\]

The claim then follows because $\cK_{X,\cM}^\bullet$ is supported on $\{x\}$.
\end{proof}

Now, by looking at the long exact sequence in cohomology associated to the triangle \eqref{eq-triangleShifted}, we get the following: 

\begin{prop} \label{prop-DescribeCoho} In our one-step resolution set-up, assume $X$ is proper. Then we have short exact sequences for any $j \geq 2$:
\[ 0 \to {\rm coker}({\rm IH}^{j-1}_{\cM}(X) \to {\rm IH}^{j-1}_{\cM,x}(X)) \to H^j_{\cM}(X) \to \ker({\rm IH}^j_{\cM}(X) \to {\rm IH}^j_{\cM,x}(X)) \to 0,\]
and an isomorphism
\[ H^1_{\cM}(X) \cong \ker({\rm IH}^1_{\cM}(X) \to {\rm IH}^1_{\cM,x}(X)).\]
\end{prop}

The left object is pure of weight $j-1$ and the right object is pure of weight $j$. Using the Lefschetz structure on $H^j_{\cM}(\widetilde{Z})$, note that the previous lemma says that the outer terms are determined by canonical morphisms
\[ {\rm IH}_{\cM}^j(X) \to H^j_{\cM,{\rm prim}}(\widetilde{Z}).\]

\section{Applications to Singularity Invariants}
The results of the previous section describe the relevant Hodge modules themselves in terms of primitive cohomology. We now translate those structural formulas into concrete statements about the singularity invariants considered earlier.

\subsection{${\rm HRH}(X)$ and $c(X)$}
We keep the same set-up and notation as in \corollaryref{cor-HodgeFiltDualObject} above, and in particular, we work with mixed Hodge modules instead of general mixed sheaves.

\begin{cor} \label{cor-kCCI} Let $k = \min\{ {\rm HRH}(\widetilde{X}), {\rm HRH}(\widetilde{Z})\}$ which we assume is non-negative (and possibly $+\infty$).

If $k \geq 1$, $\delta > 0$, and $(d,\delta) \neq (1,1)$, then
$c(X) \leq 0$, with equality if and only if
\[ F^{d-b}{\rm IH}^{d-b}_{\rm prim}(\widetilde Z)=0
\qquad \text{ for all } 0\leq b\leq d-1.\]

If $(d,\delta) = (1,1)$, then $c(X) \geq 0$ if and only if ${\rm IH}^1(\widetilde{Z}) = 0$.

If $\delta = 0$, then for any $0 \leq \ell \leq k$ we have $c(X) \geq \ell$ if and only if 
\[ F^{d-b-\ell} {\rm IH}^{d-b}_{\rm prim}(\widetilde{Z}) = 0 \text{ for all } 0 < b < d.\]
\end{cor}
\begin{proof}
For ease of notation, we write
\[ n=\dim X=d+\delta+1, \qquad P^q={\rm IH}^q_{\rm prim}(\widetilde Z).\]

Fix an integer $0\leq\ell\leq k$. Applying
\corollaryref{cor-HodgeFiltDualObject} with $k$ there replaced by
$\ell$, we see that the inequality $c(X)\geq\ell$ is equivalent to the following two
families of vanishings:
\begin{equation} \label{eq-FirstVanishing} F^{n-\ell-r-j-1}P^{d-(j-\delta+2r)}=0\quad \text{ for all } 0<j\leq n-2,\delta\leq j,\ 0\leq r\leq\delta, \end{equation}
and
\begin{equation} \label{eq-SecondVanishing}F^{n-\ell-r-\delta-1}P^{d-(\delta-j+2r)}=0 \quad \text{ for all }  0<j<\delta,\ 0\leq r\leq j. \end{equation}

Suppose first that $\delta>0$ and $d>1$. In \eqref{eq-FirstVanishing}, set
\[
b=j-\delta+2r.
\]
The corresponding condition becomes
\[ F^{d-b-\ell+r}P^{d-b}=0. \]
For fixed $b$, the strongest vanishing condition is obtained by taking $r$ as
small as possible. If $0\leq b\leq d-1$, we take
$(j,r)=(\delta+b,0)$ and obtain
\begin{equation} \label{eq-StrongestFirstVanishing}
F^{d-b-\ell}P^{d-b}=0. \end{equation} 
For $b=d$, the smallest possible value of $r$ is $1$, obtained from
$(j,r)=(\delta+d-2,1)$, and the resulting condition is
\begin{equation} \label{eq-StrongestFirstVanishingbd}
F^{1-\ell}P^0=0. \end{equation} 
Thus \eqref{eq-FirstVanishing} is equivalent to \eqref{eq-StrongestFirstVanishing}, for $0\leq b\leq d-1$, together
with \eqref{eq-StrongestFirstVanishingbd}.

The family \eqref{eq-SecondVanishing} adds no further conditions. Indeed, set
\[
b=\delta-j+2r.
\]
Its corresponding condition is
\[
F^{d-\ell-r}P^{d-b}=0.
\]
Since
\[
b-r=\delta-j+r>0,
\]
if $b\leq d-1$, this condition follows from \eqref{eq-StrongestFirstVanishing}, because
\[
d-\ell-r=(d-b-\ell)+(b-r).
\]
If $b=d$, it follows from \eqref{eq-StrongestFirstVanishingbd} and if $b>d$, the cohomological degree
is negative.

Consequently, we have the equivalence
\begin{equation} \label{eq-EquivalenceForC}
c(X)\geq\ell
\quad\iff\quad \begin{cases}
F^{d-b-\ell}P^{d-b}=0 &  \text{ for all } 0\leq b\leq d-1,\\
F^{1-\ell}P^0=0.
\end{cases} \end{equation}
Now $P^0=H^0(\widetilde Z)$ satisfies, for $\ell \geq 0$, the equivalence
\[
F^{1-\ell}P^0 = 0 \quad\iff\quad \ell=0.\]
Taking $\ell = 0$ in \eqref{eq-EquivalenceForC} gives the asserted criterion for
$c(X) \geq 0$, while taking $\ell = 1$ shows that $c(X) < 1$. Thus, when
$k \geq 1$, we have $c(X) \leq 0$, with equality precisely under the
vanishings in the statement.

Now suppose $d=1$ and $\delta > 0$. The only potentially nonzero term in \eqref{eq-FirstVanishing} is
\[ F^{1-\ell}P^1=0.\]
If $\delta>1$, the only potentially nonzero term in \eqref{eq-SecondVanishing} is obtained from
$j=\delta-1$ and $r=0$, and is
\[ F^{1-\ell} P^0 = 0.\]
The preceding argument therefore gives the first assertion in this
case as well.

If $\delta=d=1$, the family \eqref{eq-SecondVanishing} is empty, and for $\ell=0$ the only
condition is
\[
F^1{\rm IH}^1(\widetilde Z)=0.
\]
Since ${\rm IH}^1(\widetilde Z)$ is a pure Hodge structure of weight
one, Hodge symmetry shows that this is equivalent to
${\rm IH}^1(\widetilde Z)=0$.

Finally, suppose $\delta=0$. Then \eqref{eq-SecondVanishing} is empty, while in \eqref{eq-FirstVanishing} we have
$r=0$ and $b=j$, with $1\leq b\leq d-1$. Hence
\[c(X)\geq\ell
\quad\iff\quad
F^{d-b-\ell}{\rm IH}^{d-b}_{\rm prim}(\widetilde Z)=0
\qquad \text{ for all } 0<b<d,\]
as claimed.
\end{proof}

With this in hand, we can also study ${\rm HRH}(X)$.

\begin{cor} \label{cor-HRHIsolated} In the above notation, if either \[\begin{cases} \delta > 0 \text{ and } k\geq 1 \\ (d,\delta) = (1,1) \text{ and } k \geq 0\end{cases},\] then $c(X) \geq 0$ if and only if ${\rm HRH}(X) = 0$.

If $\delta = 0$, then for any $0 \leq \ell \leq k$ we have ${\rm HRH}(X) \geq \ell$ if and only if $c(X) \geq \ell$ and moreover
\[ F^{d-\ell} {\rm IH}^{d}_{\rm prim}(\widetilde{Z}) = 0.\]
\end{cor}
\begin{proof} The variety $X$ satisfies ${\rm HRH}(X) \geq \ell$ if and only if $c(X) \geq \ell$ and, moreover, if
\[ F_{\ell-\dim X}(\cH^0 \mathbf D_X^H/{\rm IC}_X^H) = 0.\]

We can apply \corollaryref{cor-HodgeFiltDualObject} to see that this last condition is equivalent (using decreasing filtrations for the Hodge structures) to the vanishing
\[ F^{d -\ell}{\rm IH}^{d-\delta}_{\rm prim} (\widetilde{Z})= 0.\]

Now, we will handle the cases as in \corollaryref{cor-kCCI} above.

First, if we assume $\delta > 0$ and $d >1$, or $\delta >1$ and $d =1$, then $c(X) \geq 0$ is equivalent to the vanishing
\[ F^{d-b} {\rm IH}^{d-b}_{\rm prim}(\widetilde{Z}) = 0 \text{ for all } 0 \leq b \leq d-1.\]

Rewriting the conditions for ${\rm HRH}(X) \geq 0$ above, we see that we must moreover require the vanishing:
\[ F^{d}{\rm IH}^{d-\delta}_{\rm prim}(\widetilde{Z}) = 0,\]
but this is automatic for $\delta >0$. This proves the equivalence in this case.

If $\delta = d =1$, then note that $c(X) \geq 0$ if and only if 
\[ {\rm IH}^1(\widetilde{Z}) = 0,\]
in which case $c(X) \geq k$. Assuming this, ${\rm HRH}(X) \geq \ell$ if and only if we additionally have $F^{1-\ell} {\rm IH}^0(\widetilde{Z}) =0$. Again, the latter is true if and only if $\ell=0$, so we get the upper bound ${\rm HRH}(X) \leq 0$. 

Finally, for $\delta = 0$, we saw above ${\rm HRH}(X) \geq \ell$ if and only if $c(X) \geq \ell$ and
\[ F^{d-\ell} {\rm IH}^{d}_{\rm prim}(\widetilde{Z}) = 0.\]
\end{proof}

\subsection{Small Resolutions}
Above we saw that if $\delta > 0$, then ${\rm HRH}(X) \leq 0$. In particular, $X$ is not a rational homology manifold. The following shows that this is not special to the one-step resolution setting, but really intrinsic to any variety admitting a small resolution. The authors would like to thank Bhargav Bhatt, who showed them the following argument:
\begin{prop} Let $(X,x)$ be an isolated singularity which admits a projective birational morphism $f \colon (\widetilde{X},E) \to (X,x)$ such that $\widetilde{X}$ is a rational homology manifold and $E$ is a positive-dimensional variety with codimension $c  >1$ in $\widetilde{X}$.

Then ${\rm HRH}(X) \leq 0$. In particular, $X$ is not a rational homology manifold.
\end{prop}
\begin{proof} We use Saito's decomposition theorem for the pure Hodge module $\Q_{\widetilde{X}}^H[n]$ which says that 
\[ f_* \Q_{\widetilde{X}}^H[n] \cong {\rm IC}_X^H \oplus i_*\Sigma,\]
where, as above, $i$ is the inclusion of $\{x\}$. By the same reasoning as above, $\cH^j(\Sigma) = 0$ for $j < \min\{0,-(n-2c)\}$.

By applying $i^*$, we get
\[ a_* \Q_E^H[n] \cong i^*({\rm IC}_X^H) \oplus \Sigma,\]
and so by applying ${\rm Gr}^F_{-1} \cH^{2-n}$, we get, from the inequality $2-n < 2c-n$, an isomorphism
\[ {\rm Gr}^F_{-1} H^{2}(E) \cong {\rm Gr}^F_{-1} \cH^{2-n} i^*({\rm IC}_X^H).\]

We have the exact triangle
\[ i_* \widetilde{\cK}_X^\bullet = \cK_X^\bullet \to \Q_X^H[n] \to {\rm IC}_X^H \xrightarrow[]{+1},\]
where $\widetilde{\cK}_X^\bullet \in D^b({\rm MHM}(\{x\})) = D^b({\rm MHS})$, and by definition, the condition ${\rm HRH}(X) \geq k$ is equivalent to ${\rm Gr}^F_{-p}{\rm DR}(\cK_X^\bullet) =0$ for all $p\leq k$. Hence, as ${\rm DR}(-)$ commutes with $i_*$, we see that ${\rm HRH}(X) \geq k$ is equivalent to
\[ {\rm Gr}^F_{-p} \widetilde{\cK}_X^\bullet = 0 \text{ for all } p \leq k.\]

So if we assume ${\rm HRH}(X) \geq 1$, then by applying ${\rm Gr}^F_{-1}\cH^{2-n}i^*(-)$ to the triangle above, we get the isomorphism
\[ {\rm Gr}^F_{-1} \cH^{2-n}i^*(\Q_X^H[n]) \cong {\rm Gr}^F_{-1}\cH^{2-n}i^*({\rm IC}_X^H),\]
but the space on the left is $0$. Indeed, $i^*(\Q_X^H[n]) = \Q^H[n]$, hence its $(2-n)$th cohomology vanishes. In conclusion, we see that ${\rm HRH}(X) \geq 1$ implies 
\[ {\rm Gr}_F^1 H^2(E) = 0,\]
but as $E$ is projective this space is never zero.
\end{proof}

\subsection{Level of Generation of Hodge filtration}
Now, assume that $X$ is embedded into a smooth variety $Y$ via $\iota \colon X \to Y$, which is a closed embedding of codimension $q = \dim Y - n$. Assume that $\widetilde{X}$ and $\widetilde{Z}$ are rational homology manifolds. 

The local cohomology mixed Hodge modules are computed by \lemmaref{lem-DXLocCoh} above.

As the $c(X)$ invariant measures where the Hodge filtration on local cohomology starts, it is also natural to ask where it ``ends'', to be made precise by the \emph{generation level}.

If $(\cM,F)$ is a filtered left $\cD_Y$-module underlying a mixed Hodge module on the smooth variety $Y$, then by definition, $F_\bullet \cM$ is a good filtration. In particular, there exists some $k_0$ such that, for all $\ell \geq 0$, we have equality
\[ F_\ell \cD_Y \cdot F_{k_0}\cM = F_{k_0+\ell}\cM.\]

If this equality holds, we say $F_\bullet \cM$ is generated at level $k_0$. 

By definition of the filtered de Rham complex, we get the following:
\begin{lem}[\cite{MPLocCoh}*{Lem. 10.1}] \label{lem-genLevelGrDR} If $(\cM,F)$ is a filtered left $\cD_Y$-module underlying a mixed Hodge module on the smooth variety $Y$, then $F_\bullet \cM$ is generated at level $k_0$ if and only if $\cH^0 {\rm Gr}^F_{i-\dim Y} {\rm DR}_Y(\cM) = 0$ for all $i > k_0$.
\end{lem}

\begin{rmk} The generating levels of the Hodge filtration on local cohomology modules as well as that on nearby/vanishing cycles (in the hypersurface case) are related in \cites{MPLocCoh,CDM,OR} to various singularity invariants.
\end{rmk}

The following definition differs from that in \cite{SaitoOnTheHodge} by a shift of $p(\cM) = \min\{p\mid F_p \cM \neq 0\}$.

\begin{defi} We say that $(\cM,F)$ has \emph{generation level} $k$ if $F_\bullet \cM$ is generated at level $k$ but not $k-1$. We denote the generation level by ${\rm gl}(\cM,F)$.
\end{defi}

\begin{lem} \label{lem-genLevelGrDRPrecise} Let $(\cM,F)$ be a non-zero filtered left $\cD_Y$-module underlying a mixed Hodge module on $Y$. Then
\[ {\rm gl}(\cM,F) = \max\{p \mid \cH^0 {\rm Gr}^F_{p-\dim Y} {\rm DR}(\cM) \neq 0\}.\]
\end{lem}

We would like to give a bound on the index at which the Hodge filtrations of \[\cH^{q+j}_X(\cO_Y) = \iota_* \cH^j(\mathbf D_X^H)(-q) \text{ and } {\rm Gr}^W_{n+2q+1}\cH^{q}_X(\cO_Y) =\iota_* {\rm Gr}^W_{n+1}\cH^0(\mathbf D_X^H)(-q)\] are generated. The bound we get resembles that in \cite{MPLocCoh}*{Thm. 10.2}, though their result is only for the right-most non-vanishing local cohomology (and holds in general).

Let $\mu^\ell_{\rm prim}(\widetilde{Z}) = \sup\{ p \mid F^p H^\ell_{\rm prim}(\widetilde{Z}) = H^{\ell}_{\rm prim}(\widetilde{Z})\}$, which satisfies the trivial inequality $\mu^{\ell}_{\rm prim}(\widetilde{Z}) \geq 0$.

\begin{cor} \label{cor-GenLevel} In the above set-up, we have 
\[{\rm gl}({\rm Gr}^W_{n+2q+1} \cH^q_X(\cO_Y),F) = d - \mu^{d-\delta}_{\rm prim}(\widetilde{Z})\].

For $0 < j \leq \dim X -2$, we have
\[ {\rm gl}(\cH^{q+j}_X(\cO_Y),F) = \begin{cases} d - (j-\delta) - \min_{0 \leq r \leq \delta}\{ \mu_{\rm prim}^{d-(j-\delta+2r)}(\widetilde{Z}) + r\} & \delta \leq j \\ d - \min_{0 \leq r \leq j}\{ \mu_{\rm prim}^{d - (\delta-j+2r)}(\widetilde{Z})+ r\} & \delta > j \end{cases}.\]
\end{cor}
\begin{proof} Let $i \colon Z \to Y$ be the closed embedding. We have the isomorphisms for $j>0$
\[ \cH^{q+j}_X(\cO_Y) \cong \begin{cases} \bigoplus_{r=0}^{\delta} i_* H^{d-(j-\delta+2r)}_{\rm prim}(\widetilde{Z})(-q-j-r-1) & \delta \leq j \\ \bigoplus_{r = 0}^{j} i_*H^{d-(\delta-j+2r)}_{\rm prim}(\widetilde{Z})(-q-\delta-r-1) & \delta >j \end{cases},\]
and for $j=0$, we have
\[{\rm Gr}^W_{n+2q+1} \cH^q_X(\cO_Y) \cong i_* H^{d-\delta}_{\rm prim}(\widetilde{Z})(-q-\delta-1).\]

Now, by applying $\cH^0 {\rm Gr}^F_{p- \dim Y} {\rm DR}(-)$ and using the commutativity with $i_*$, we get (writing decreasing Hodge filtrations on the Hodge structure on the right hand side):
\[ \cH^0 {\rm Gr}^F_{p- \dim Y} {\rm DR}(\cH^{q+j}_X(\cO_Y)) \cong \begin{cases} \bigoplus_{r=0}^{\delta} i_* {\rm Gr}_F^{\dim Y -p-q-j-r-1}(H^{d-(j-\delta+2r)}_{\rm prim}(\widetilde{Z})) & \delta \leq j \\ \bigoplus_{r = 0}^{j} i_*{\rm Gr}_F^{\dim Y - p-q-\delta-r-1}(H^{d-(\delta-j+2r)}_{\rm prim}(\widetilde{Z})) & \delta >j \end{cases},\]
and
\[\cH^0 {\rm Gr}^F_{p- \dim Y} {\rm DR}({\rm Gr}^W_{n+2q+1} \cH^q_X(\cO_Y)) \cong i_* {\rm Gr}_F^{\dim Y - p-q-\delta-1}(H^{d-\delta}_{\rm prim}(\widetilde{Z})).\]

The claim now follows by \lemmaref{lem-genLevelGrDRPrecise}.
\end{proof}

For the same embedding $X\subseteq Y$, we study the generating level of $\iota_*{\rm IC}_X^H(-q) = W_{n+2q} \cH^q_X(\cO_Y)$. We assume moreover in this section that $\widetilde{X}$ is smooth.

In the following diagrams and formulas, we suppress the closed-embedding pushforward $i_*$ on objects supported on $Z$.

Recall that by the decomposition theorem, we have an isomorphism
\[ \iota_* f_*\Q^{H}_{\widetilde{X}}[\dim X] \cong \iota_* {\rm IC}_X^H \oplus  i_* H^{\dim X}(\widetilde{Z}) \oplus \bigoplus_{\ell =1}^{d-c_{\widetilde{Z}}} \left(i_*H^{\dim X +\ell}(\widetilde{Z})[-\ell] \oplus i_*H^{\dim X+\ell}(\widetilde{Z})(\ell)[\ell]\right).\]

If we push forward to $Y$ and then apply $\cH^0 {\rm Gr}^F_{i-\dim Y} {\rm DR}(-)$ to both sides, we get
\[ \iota_* \cH^0 Rf_* \Omega_{\widetilde{X}}^{\dim Y -i}[\dim X +i - \dim Y] \cong \cH^0 {\rm Gr}^F_{i-\dim Y} {\rm DR}({\rm IC}_X^H) \oplus {\rm Gr}_F^{\dim Y - i} H^{\dim X}(\widetilde{Z}),\]
where we used that ${\rm Gr}^F_p{\rm DR}(-)$ on a point agrees with ${\rm Gr}_F^{-p}(-)$ and only has cohomology in degree $0$. This can be rewritten as
\[ \iota_*R^{i-q}f_* \Omega_{\widetilde{X}}^{\dim X -(i-q)} \cong \cH^0 {\rm Gr}^F_{i-\dim Y} {\rm DR}({\rm IC}_X^H) \oplus {\rm Gr}_F^{\dim X - (i-q)} H^{\dim X}(\widetilde{Z}),\]
and if we instead study ${\rm IC}_X^H(-q)$, then we can replace the $i$ by $i+q$ in this formula and get the cleaner isomorphism
\[ \iota_*R^{i}f_* \Omega_{\widetilde{X}}^{\dim X - i} \cong \cH^0 {\rm Gr}^F_{i-\dim Y} {\rm DR}({\rm IC}_X^H(-q)) \oplus {\rm Gr}_F^{\dim X - i} H^{\dim X}(\widetilde{Z}).\]

By the Lefschetz isomorphisms, we have
\[ {\rm Gr}_F^{\dim X - i} H^{\dim X}(\widetilde{Z}) \cong {\rm Gr}_F^{d-i} H^{d-c_{\widetilde{Z}}}(\widetilde{Z}),\]
and so the second term vanishes automatically for $i > d$. 

We conclude the following:
\begin{cor} In the above notation, we have the inequality
\[ d  \geq \max \{ k \mid R^k f_* \Omega_{\widetilde{X}}^{\dim X-k} \neq 0\} \geq {\rm gl}({\rm IC}_X^H(-q),F)\]
and thus,
\[ d \geq {\rm gl}(\cH^q_X(\cO_Y),F).\]
\end{cor}
\begin{proof} The first inequality follows from the fact that $\widetilde{Z}$ is the maximal dimensional fiber of the map $f$, and it has dimension $d$ \cite{HartshorneBook}*{Cor. III.11.2}. The second inequality follows from \lemmaref{lem-genLevelGrDRPrecise}.

The last line follows from the short exact sequence
\[ 0 \to {\rm IC}_X^H(-q) \to \cH^q_X(\cO_Y) \to {\rm Gr}^W_{n+2q+1} \cH^q_X(\cO_Y) \to 0,\]
inducing the exact sequence
\[ \cH^0 {\rm Gr}^F_{i-\dim Y} {\rm DR}({\rm IC}_X^H(-q)) \to \cH^0 {\rm Gr}^F_{i-\dim Y} {\rm DR}(\cH^q_X(\cO_Y)) \to \cH^0 {\rm Gr}^F_{i-\dim Y} {\rm DR}({\rm Gr}^W_{n+2q+1} \cH^q_X(\cO_Y)).\]

The left-most term vanishes for $i > d$ by the first claim of this corollary and for the right-most term it is $d- \mu^{d-\delta}_{\rm prim}(\widetilde{Z})$ by the first claim of \corollaryref{cor-GenLevel}. By non-negativity of $\mu^{d-\delta}_{\rm prim}(\widetilde{Z})$, this proves the claim.
\end{proof}

\subsection{Hodge--Lyubeznik Numbers}
We now provide the application to Hodge--Lyubeznik numbers. We will use this below to recover the computation of \cite{HodgeLyubeznik}*{Example (1)} and extend it to cones over projective rational homology manifolds.

Recall that we are assuming $Z = \{x\}$ is a point, $\widetilde{Z}$ and $\widetilde{X}$ are rational homology manifolds, and $n\geq 2$.

We define primitive Hodge numbers for $\widetilde{Z}$ for $0 \leq k \leq d$ by
\[ \underline{h}^{p,k-p}_{\rm prim}(\widetilde{Z}) = \dim_\C {\rm Gr}^F_{-p} H^k(\widetilde{Z})_{\rm prim},\]
which is defined to mimic the standard notation
\[ \underline{h}^{p,k-p}(\widetilde{Z}) = \dim_{\C} {\rm Gr}^F_{-p} H^k(\widetilde{Z}).\]

For $y\neq x$, we can apply \lemmaref{lem-HLRHM} above to compute $\lambda_{r,s}^{p,q}(\cO_{X,y})$, using the fact that $X$ is CCI away from $x$.

The goal of this subsection is to prove the following theorem:

\begin{thm} \label{thm-HLIsolated} For $2\leq s < n$, the Hodge--Lyubeznik number $\lambda_{r,s}^{p,q}(\cO_{X,x})$ is non-zero only if $r= 0$ and $p+q=1-s$, in which case it is equal to
\[ \lambda_{0,s}^{p,q}(\cO_{X,x}) = \begin{cases} \sum_{a=0}^{\delta} \underline{h}^{-q-a,-p-a}_{\rm prim}(\widetilde{Z}) & \delta \leq n -s \\ \sum_{a = 0}^{n-s} \underline{h}^{n-s-\delta-q-a,n-s-\delta-p-a}_{\rm prim}(\widetilde{Z}) & \delta > n -s \end{cases}.\]

For $s = n$, we have $\lambda_{r,n}^{p,q}(\cO_{X,x}) = 0$ for $r\leq 1$. For $r\geq 2$, we have
\[ \lambda_{r,n}^{p,q}(\cO_{X,x}) = {\rm I}\lambda_r^{p,q}(\cO_{X,x}).\]

The intersection Hodge--Lyubeznik numbers are non-zero only for $r\geq 1$ and $p+q = r-n$, and are given by the following:
\[ {\rm I}\lambda_r^{p,q}(\cO_{X,x}) = \begin{cases} \sum_{a=0}^{\delta} \underline{h}^{-q-a,-p-a}_{\rm prim}(\widetilde{Z}) & r > \delta \\ \sum_{a=0}^{r-1} \underline{h}^{d-a+p,d-a+q}_{\rm prim}(\widetilde{Z}) & r \leq \delta \end{cases}.\]
\end{thm}
\begin{proof}
First of all, for $2 \leq s < n$ (the modules are $0$ for $s =0, 1$ by the bound on the local cohomological defect), we have that $\cH^{n-s} \mathbf D_X^H(n) = i_{x,*} V_s$ is supported on $\{x\}$, hence we can write
\[ i_x^! \cH^{n-s} \mathbf D_X^H(n) = V_s,\]
and in particular, we see that $\cH^r_x \cH^{n-s} \mathbf D_X^H(n) \neq 0$ implies $r = 0$. 

Moreover, we see that
\[ {\rm Gr}^W_{p+q} \cH^0_x \cH^{n -s} \mathbf D_X^H(n) = {\rm Gr}^W_{p+q} V_s = \cH^0_x {\rm Gr}^W_{p+q} \cH^{n-s} \mathbf D_X^H(n),\]
and so we can immediately use \corollaryref{cor-Isolated}. We will use $j$ in place of $n -s$ for ease of notation below. 

For $r=0$, by Tate twisting that result by $n$, we get for $0< j \leq n-2$ the isomorphisms
\[ {\rm Gr}^W_{-n+j+1} \cH^0_x \cH^{j} \mathbf D_X^H(n) \cong \begin{cases} \bigoplus_{a=0}^{\delta} H^{d-(j-\delta+2a)}_{\rm prim}(\widetilde{Z})(n-a-1-j) & \delta \leq j \\ \bigoplus_{a = 0}^{j} H^{d-(\delta-j+2a)}_{\rm prim}(\widetilde{Z})(d-a) & \delta > j \end{cases},\]
and for such $j$ this is the only non-zero ${\rm Gr}^W$ piece.

If we apply ${\rm Gr}^F_{-p}$, then we get
\[ {\rm Gr}^F_{-p}{\rm Gr}^W_{-n+j+1} \cH^0_x \cH^{j} \mathbf D_X^H(n)\cong \begin{cases} \bigoplus_{a=0}^{\delta} {\rm Gr}^F_{a+j+1-n-p} H^{d-(j-\delta+2a)}_{\rm prim}(\widetilde{Z}) & \delta \leq j \\ \bigoplus_{a = 0}^{j} {\rm Gr}^F_{a-d-p}H^{d-(\delta-j+2a)}_{\rm prim}(\widetilde{Z}) & \delta > j \end{cases}.\]

If we want to express this in terms of the primitive Hodge numbers, we negate the index of the Hodge filtration, and then we subtract the result from the weight. Thus, the second term in the primitive Hodge number is
\[ \begin{cases} (d+\delta-j-2a)+ (a+j+1-n -p) = -p-a & \delta \leq j \\ (d+j-\delta-2a) + (a-d-p) = j-\delta-p-a & \delta > j\end{cases}.\]

We thus get equality of $\dim_{\C} {\rm Gr}^F_{-p} {\rm Gr}^W_{p+q} \cH^0_x \cH^j \mathbf D_X^H(n)$ (by setting $p+q = -n+j+1$) with
\[ \begin{cases} \sum_{a=0}^{\delta} \underline{h}^{n+p-a-j-1,-p-a}_{\rm prim}(\widetilde{Z})& \delta \leq j \\ \sum_{a = 0}^{j} \underline{h}^{d+p-a,j-\delta-a-p}_{\rm prim}(\widetilde{Z}) & \delta > j \end{cases}.\]

Finally, letting $j = n -s$ for $2\leq s < n$, we get that the Hodge--Lyubeznik number is non-zero only if $r = 0$ and $p+q = 1-s$, in which case it is equal to
\[ \lambda_{0,s}^{p,q}(\cO_{X,x}) = \begin{cases} \sum_{a=0}^{\delta} \underline{h}^{-q-a,-p-a}_{\rm prim}(\widetilde{Z}) & \delta \leq n-s \\ \sum_{a = 0}^{n-s} \underline{h}^{n-s-\delta-q-a,n-s-\delta-p-a}_{\rm prim}(\widetilde{Z}) & \delta > n-s \end{cases}.\]

For $s= n$, our goal is to compute
\[ \lambda_{r,n}^{p,q}(\cO_{X,x}) = \dim_{\C} {\rm Gr}^F_{- p} {\rm Gr}^W_{p+q} \cH^r_x(\cH^{0} \mathbf D_X^H(n)) \]
which we can do in two steps. Indeed, we have the short exact sequence
\[ 0 \to {\rm IC}_X^H \to \cH^0 \mathbf D_X^H \to \cH^0 \mathbf D_X^H/{\rm IC}_X^H \to 0,\]
to which we apply $i_x^!$. As noted above, $\cH^0 \mathbf D_X^H/{\rm IC}_X^H = i_{x,*} V$ is supported on $Z = \{x\}$ when we assume $\widetilde{X}$ is a rational homology manifold. Thus, when we apply $i_x^!$, we get isomorphisms
\begin{equation} \label{eq-isoHL} \cH^j_x {\rm IC}_X^H(n) \cong \cH^j_x \cH^0\mathbf D_X^H(n) \text{ for all } j > 1.\end{equation}

As far as $\mathbf D_X^H$ is concerned, we see that those are the only interesting pieces:

\begin{lem} We have $\cH^0_x {\rm IC}_X^H = \cH^0_x \cH^0 \mathbf D_X^H = \cH^1_x \cH^0 \mathbf D_X^H = 0$.
\end{lem}
\begin{proof} The vanishing $\cH^0_x {\rm IC}_X^H = \cH^0_x \cH^0 \mathbf D_X^H = 0$ is true because both ${\rm IC}_X^H$ and $\cH^0\mathbf D_X^H$ have no sub-objects supported on $\{x\}$. To see this is true for $\cH^0 \mathbf D_X^H$, note that for $M \in {\rm MHS}$, we have by the usual $t$-structure property:
\[ {\rm Hom}(i_* M, \cH^0\mathbf D_X^H)= {\rm Hom}(i_* M, \mathbf D_X^H),\]
then, by adjunction, this is equal to
\[ {\rm Hom}(M,i^! \mathbf D_X^H) = {\rm Hom}(M,\Q^H(-\dim X)[-\dim X]) = {\rm Ext}^{-\dim X}(M, \Q^H(-\dim X))= 0,\]
where we use the isomorphism \eqref{eq-DXRule} and the fact that the Ext group between objects in an abelian category with negative index vanishes.

To see the vanishing $\cH^1_x \cH^0 \mathbf D_X^H = 0$, consider the spectral sequence
\[ E_2^{p,q} = \cH^p_x \cH^q \mathbf D_X^H \implies \cH^{p+q}_x \mathbf D_X^H.\]

By the isomorphism \eqref{eq-DXRule}, we see that $\cH^{p+q}_x \mathbf D_X^H = 0$ unless $p+q = n$. Hence $E_\infty^{p,q} = 0$ for all $p+q \neq n$.

On the other hand, using that for $q>0$ the module $\cH^q \mathbf D_X^H$ is supported on $Z=\{x\}$, we have $E_2^{p,q} =0 $ if $q >0$ and $p\neq 0$. If $q = 0$, then $E_2^{p,q} = 0$ unless $p>0$ by the first vanishing that we argued. It is easy to see then that $E_2^{1,0} = E_\infty^{1,0}$.

As we are assuming $c_Z \geq 2$, we have that $n > 1$, and so we get the vanishing $\cH^{1}_x \cH^0\mathbf D_X^H = E_2^{1,0} = E_\infty^{1,0} = 0$, as claimed.
\end{proof}

We have, by \corollaryref{cor-Isolated}, the decomposition:
\[ f_* \Q_{\widetilde{X}}^H[n] \cong {\rm IC}_X^H \oplus  i_*H^n(\widetilde{Z}) \oplus \bigoplus_{\ell =1}^{d-c_{\widetilde{Z}}} \left(i_*H^{n+\ell}(\widetilde{Z})[-\ell] \oplus i_*H^{n+\ell}(\widetilde{Z})(\ell)[\ell]\right).\]
Combined with base change we have, by applying $i_x^*$, the decomposition:
\[a_*\Q_{\widetilde{Z}}^H[n] \simeq i_x^*{\rm IC}_X^H\oplus H^n(\widetilde{Z}) \oplus \bigoplus_{\ell=1}^{d-c_{\widetilde{Z}}} \left(H^{n+\ell}(\widetilde{Z})[-\ell] \oplus H^{n+\ell}(\widetilde{Z})(\ell)[\ell] \right).\]

Taking cohomology we get
\[ H^{n-j}(\widetilde{Z}) \simeq \cH^{-j} i_x^*{\rm IC}_X^H \oplus \begin{cases} H^{n-j}(\widetilde{Z}) & j \leq   0 \\ H^{n+j}(\widetilde{Z})(j) & j > 0 \end{cases}.\]

We see $\cH^{-j} i_x^*{\rm IC}_X^H = 0$ for $j\leq 0$, but we already knew that because ${\rm IC}_X^H$ has no quotient objects supported on $\{x\}$.

We consider $j>0$, but we are actually interested in $i_x^! {\rm IC}_X^H$, and so we must dualize. By pure polarizability, this just amounts to tracking some Tate twists. We get for $j>0$
\begin{equation} \label{eq-IsoBeforeLefschetz} H^{n-j}(\widetilde{Z})(n-j) \cong \cH^j i_x^! {\rm IC}_X^H(n) \oplus H^{n+j}(\widetilde{Z})(n).\end{equation}

By the Lefschetz structure on the cohomology of $\widetilde{Z}$, we can simplify the isomorphism \eqref{eq-IsoBeforeLefschetz}. Indeed, we have the surjection
\[ c_1(L\vert_{\widetilde{Z}})^j \colon H^{n-j}(\widetilde{Z})(n-j) \to H^{n+j}(\widetilde{Z})(n), \]
and so for $j>0$, we can identify $\cH^j i_x^! {\rm IC}_X^H(n)$ with the kernel of this map. We can rewrite $n-j = d +\delta -(j-1)$, so that we have the Lefschetz decomposition
\[H^{n-j}(\widetilde{Z}) \cong \begin{cases} \bigoplus_{a\geq 0} H^{d-(j-1-\delta+2a)}_{\rm prim}(\widetilde{Z})(-a) & \delta \leq j-1\\ \bigoplus_{a\geq 0}H^{d-(\delta - (j-1)+2a)}_{\rm prim}(\widetilde{Z})(j-a-\delta-1) & \delta > j-1\end{cases},\]
and as before we can use this to compute the kernel of $c_1(L\vert_{\widetilde{Z}})^j$.

In summary (rewriting $r$ in place of $j$), we have for $r>0$ the isomorphisms
\[ {\rm Gr}^W_{r-n} \cH^r i_x^! {\rm IC}_X^H (n) \cong  \begin{cases}\bigoplus_{a= 0}^{\delta} H^{d-(r-1-\delta+2a)}_{\rm prim}(\widetilde{Z})(n-r-a) & r > \delta \\ \bigoplus_{a= 0}^{r-1} H^{d-(\delta - r +2a+1)}_{\rm prim}(\widetilde{Z})(n-a-\delta-1) & r\leq \delta\end{cases}\]
and all other ${\rm Gr}^W$ pieces are $0$.

We can use these results to compute the (intersection) Hodge--Lyubeznik numbers in terms of the primitive Hodge numbers of $\widetilde{Z}$, similarly to what we did for $s < n$, though we omit the computation as it is the same as above. 
\end{proof}

\begin{rmk} With the computation of the Hodge--Lyubeznik numbers, one can reprove \corollaryref{cor-kCCI} above, using \theoremref{thm-MaincX}\eqref{thm-HLCCI}. However, this approach does not simplify that proof.
\end{rmk}

\section{Cones and Contractions of the Zero Section} \label{sec-Contractions}
We conclude by specializing the preceding results to a natural class of examples arising from contractions of zero sections. Besides producing explicit examples, this family makes the role of the parameter $\delta$ explicit, as it directly relates to the rank of the ample vector bundle.

Again, we work over $K$ a subfield of $\C$ and with a theory of mixed sheaves $\cM(-)$ on $K$-varieties.

Let $Y$ be a projective variety of positive dimension with an ample vector bundle $\cE$ of rank $e$. Below we use the definition of ample vector bundle as in \cite{Hartshorne}, though note that \emph{loc. cit.} assumes the ground field is algebraically closed. As we are mostly interested in the behavior after base-change to $\C$, we do not need to worry about this caveat.

By \cite{Hartshorne}*{Prop. 3.5}, we have the diagram
\[ \begin{tikzcd} Y \ar[d,"p"] \ar[r] & \widetilde{X} \ar[d,"f"] \\ 0 \ar[r] & X \end{tikzcd}\]
where $\widetilde{X}$ is the total space of $\cE^*$, the top horizontal morphism is the inclusion of the zero section, and $f$ is a projective birational morphism which is an isomorphism away from $0$. Here 
\[ X = {\rm Spec}(\bigoplus_{\ell \geq 0} H^0({\rm Sym}^\ell(\cE))).\]

In this situation, $\widetilde{Z} = Y$ so that $c_Z = n = \dim Y +e$ and $c_{\widetilde{Z}} = e$. Assume moreover that $Y\times_K \C$ is an irreducible rational homology manifold. As $\widetilde{X}$ is a vector bundle over $Y$, it is also an irreducible rational homology manifold after complexification. 

Note that, in the notation of \theoremref{thm-general}, $e = \delta+1$, so $\delta = 0$ corresponds to line bundles whereas higher-rank contractions fit into the $\delta > 0$ case.

We first define a relatively ample line bundle for the morphism $f \times_K \C \colon \widetilde{X}\times_K \C \to X \times_K \C$. For this part of the discussion, we drop the base-change to $\C$ from the notation. As $\widetilde{X}$ is the total space of $\cE^*$, we have the projection $\pi \colon \cE^* \to Y$, and we can use the line bundle $\pi^*{\rm det}(\cE)$. By \cite{Lazarsfeld}*{1.7.8}, to prove that this line bundle is relatively ample it suffices to prove that its restriction to all fibers is ample. The fiber over any non-zero point is simply a point, so the claim is trivially true for such fibers. The fiber over $0$ is $Y$, and the restriction of $\pi^*{\rm det}(\cE)$ to $Y$ is simply ${\rm det}(\cE)$, which by \cite{Hartshorne}*{Prop. 2.6} is ample.

Then we get the following information about the singularities of $X$ by \corollaryref{cor-Isolated}.

\begin{thm} \label{thm-contractionLocCoh} Let $n = \dim Y +e = \dim X$. 

Then
\begin{enumerate} \item ${\rm lcdef}(X) = {\rm lcdef}_{\rm gen}(X) = {\rm lcdef}_{\rm gen}^{>0}(X) \leq \dim X -2$.
    \item For all $0 < j \leq \dim X -2$, we have an isomorphism of pure objects of weight $n+j+1$:
    \[ \cH^{-j} \cK_{X, \cM}^\bullet(-j-1) \cong \cH^j \mathbf D^\cM_X \cong \begin{cases} \bigoplus_{r=0}^{e-1} i_*H_{\cM,\rm prim}^{\dim Y-(j-e+2r+1)}(Y)(-j-r-1) & e-1 \leq j \\ \bigoplus_{r = 0}^{j} i_*H_{\cM,\rm prim}^{\dim Y-(e-j+2r-1)}(Y)(-e-r) & e-1 >j \end{cases}.\]
    \item We have ${\rm Gr}^W_i \cH^0\mathbf D_X^\cM \neq 0$ implies $i \in \{n,n+1\}$, and an isomorphism
    \[ \cH^0 \cK_{X,\cM}^\bullet (-1) \cong {\rm Gr}^W_{n+1} \cH^0 \mathbf D_X^\cM \cong i_* H^{\dim Y -e+1}_{\cM,\rm prim}(Y)(-e).\]

\item We have an isomorphism in $D^b \cM(X)$ 
\[ f_* A_{\widetilde{X}}^\cM[n] \cong {\rm IC}_X^\cM \oplus  i_* H^n_{\cM}(Y) \oplus \bigoplus_{\ell =1}^{\dim Y-e} \left(i_*H^{n+\ell}_{\cM}(Y)[-\ell] \oplus i_*H^{n+\ell}_{\cM}(Y)(\ell)[\ell]\right).\]
\end{enumerate}
\end{thm}

The case $e=1$ in the following corollary was already observed in \cite{DOR1}*{Prop. 9.1}, in the case $Y$ is smooth and $\cE = L$ is a line bundle. For the next three corollaries, we take $A = \Q, K =  \C, \cM(-) = {\rm MHM}(-)$. For the first two, we do not need to assume $Y$ is a rational homology manifold.

\begin{cor} \label{cor-ContractCCI} Assume ${\rm HRH}(Y) \geq k$ and let $0 \leq \ell \leq k$.

If $k\geq 1$, $e >1$ and $(\dim Y, e) \neq (1,2)$, then $c(X) \leq 0$, with equality if and only if
\[ F^{\dim Y-b} {\rm IH}^{\dim Y-b}_{\rm prim}(Y) = 0 \text{ for all } 0 \leq b\leq \dim Y-1.\]

If $(\dim Y,e) = (1,2)$, then $c(X) \geq 0$ if and only if $X$ is CCI if and only if ${\rm IH}^1(Y) = 0$.

If $e = 1$, then we have $c(X) \geq \ell$ if and only if 
\[ F^{\dim Y -b -\ell} {\rm IH}^{\dim Y-b}_{\rm prim}(Y) = 0 \text{ for all } 0 < b < \dim Y.\]

\end{cor}
\begin{proof} This is an immediate application of \corollaryref{cor-kCCI} and \corollaryref{cor-HRHIsolated}.
\end{proof}

\begin{cor} \label{cor-ContractHRH} Assume ${\rm HRH}(Y) \geq k$ and let $0 \leq \ell \leq k$.

In the above notation, if either \[\begin{cases} e > 1 \text{ and } k\geq 1 \\ (\dim Y,e) = (1,2) \text{ and } k \geq 0\end{cases},\] then $c(X) \geq 0$ if and only if ${\rm HRH}(X) = 0$.

If $e =1$, then $X$ satisfies ${\rm HRH}(X) \geq \ell$ if and only if $c(X) \geq \ell$ and, moreover,
\[ F^{\dim Y-\ell} {\rm IH}^{\dim Y}_{\rm prim}(Y) = 0.\]
\end{cor}

We have the computation of the Hodge--Lyubeznik numbers, which follows immediately from \theoremref{thm-HLIsolated}. This recovers the result of \cite{HodgeLyubeznik}*{Ex.(1)} when $Y$ is smooth and $\cE = L$ is an ample line bundle. 

\begin{cor} Assume $Y$ is a rational homology manifold.

For $2\leq s < n$, the Hodge--Lyubeznik number $\lambda_{r,s}^{p,q}(\cO_{X,0})$ is non-zero only if $r= 0$ and $p+q=1-s$, in which case it is equal to
\[ \lambda_{0,s}^{p,q}(\cO_{X,0}) = \begin{cases} \sum_{a=0}^{\delta} \underline{h}^{-q-a,-p-a}_{\rm prim}(Y) & \delta \leq n -s \\ \sum_{a = 0}^{n-s} \underline{h}^{\dim Y+1-s-q-a,\dim Y +1-s-p-a}_{\rm prim}(Y) & \delta > n -s \end{cases}.\]

For $s = n$, we have $\lambda_{r,n}^{p,q}(\cO_{X,0}) = 0$ for $r\leq 1$. For $r\geq 2$, we have
\[ \lambda_{r,n}^{p,q}(\cO_{X,0}) = {\rm I}\lambda_r^{p,q}(\cO_{X,0}).\]

The intersection Hodge--Lyubeznik numbers are non-zero only for $r\geq 1$ and $p+q = r-n$, and are given by the following: for $p+q +n = r$, we have
\[ {\rm I}\lambda_r^{p,q}(\cO_{X,0}) = \begin{cases} \sum_{a=0}^{\delta} \underline{h}^{-q-a,-p-a}_{\rm prim}(Y) & r > \delta \\ \sum_{a=0}^{r-1} \underline{h}^{\dim Y-a+p,\dim Y-a+q}_{\rm prim}(Y) & r \leq \delta \end{cases}.\]
\end{cor}

\begin{rmk} Our computation shows that the Hodge--Lyubeznik numbers of a cone over a rational homology manifold do not depend on the chosen ample $L$, similarly to \cite{RSW}*{Cor. 1.8}. This was already shown when $Y$ is smooth in \cite{HodgeLyubeznik}.
\end{rmk}

\begin{rmk} An alternative approach could be to establish the result for $\cE = L$ an ample line bundle, and then to use the fact that for $\cE$ an ample vector bundle of higher rank, we have the diagram
\[ \begin{tikzcd} \P(\cE) \ar[r] \ar[d] & \widetilde{X}'\ar[d,"f"]\\ \{0\} \ar[r] & X\end{tikzcd},\]
where $\widetilde{X}'$ is the total space of $\cO_{\P(\cE)}(-1)$ on $\P(\cE)$.
\end{rmk}

\begin{eg} Let $X \subseteq \A^4$ be defined by $xy-zw$. This is known to have a small resolution (obtained by blowing up the non-Cartier divisor $\{x=z=0\}$) with a diagram
\[ \begin{tikzcd} \P^1 \ar[r] \ar[d] & \widetilde{X} \ar[d,"f"]\\ \{0\} \ar[r] & X \end{tikzcd},\]
where $\widetilde{X}$ is smooth and $f$ is an isomorphism away from $0$. In fact, this situation fits into the above construction: if we take $\cE = \cO(1)\oplus \cO(1)$ on $\P^1$, then this is an ample vector bundle of rank $2$ on $\P^1$ such that $X$ is isomorphic to the cone over $\P(\cE)$ with conormal bundle $\cO_{\P(\cE)}(1)$.

As $X$ has hypersurface singularities, it is CCI, hence $c(X) =+\infty$. By \corollaryref{cor-ContractHRH} we have that ${\rm HRH}(X) = 0$, but $X$ is not a rational homology manifold. We can see this in another way: as $X$ has hypersurface singularities, if $b_f(s)$ is the Bernstein-Sato polynomial for $f= xy-zw$, then it is easy to see (for example, by changing coordinates so that $X$ is defined by $x_1^2+x_2^2+x_3^2+x_4^2$) that $b_f(s) = (s+1)(s+2)$. As $f$ has homogeneous isolated singularities, its spectrum is the same as its negated (reduced) Bernstein-Sato roots. Hence, we have (in the notation of \cite{DOR2}) the equality $\widetilde{\alpha}_{\Z}(f) = {\rm Sp}_{\Z,\min}(f) = 2$. Thus, by \cite{DOR2}*{Thm. A}, we conclude
\[ {\rm HRH}(X) = {\rm Sp}_{\Z,\min}(f) -2 = 0.\]
\end{eg}

\bibliography{bib}

\begin{bibdiv}
\begin{biblist}

\bib{Burke}{article}{
      author={Burke, Andrew},
       title={Local {C}ohomological {D}efect and a {C}onjecture of
  {M}ustata-{P}opa},
        date={2026},
     journal={preprint arXiv: 2602.05197},
}

\bib{CDM}{article}{
      author={Chen, Qianyu},
      author={Dirks, Bradley},
      author={{Musta\c{t}\u{a}}, Mircea},
       title={The minimal exponent and {$k$}-rationality for local complete
  intersections},
        date={2024},
     journal={J. {\'{E}}c. polytech. {Math.}},
      volume={11},
       pages={849\ndash 873},
}

\bib{CDORSecant}{article}{
      author={Chen, Qianyu},
      author={Dirks, Bradley},
      author={Olano, Sebasti\'{a}n},
      author={Raychaudhury, Debaditya},
       title={Hodge theory of secant varieties},
        date={2025},
     journal={In preparation},
}

\bib{DOR1}{article}{
      author={Dirks, Bradley},
      author={Olano, {Sebasti\'{a}n}},
      author={Raychaudhury, Debaditya},
       title={A {Hodge} {Theoretic} generalization of {$\mathbb Q$}-{Homology}
  {Manifolds} {I}: {G}eneral {C}ase},
        date={2025},
     journal={preprint arXiv:2501.14065},
}

\bib{DOR2}{article}{
      author={Dirks, Bradley},
      author={Olano, {Sebasti\'{a}n}},
      author={Raychaudhury, Debaditya},
       title={A {Hodge} {Theoretic} generalization of {$\mathbb Q$}-{Homology}
  {Manifolds} {II}: {L}ocal {C}omplete {I}ntersections},
        date={2025},
     journal={preprint arXiv:},
}

\bib{DaoTakagi}{article}{
      author={Dao, Hailong},
      author={Takagi, Shunsuke},
       title={On the relationship between depth and cohomological dimension},
        date={2016},
     journal={Compositio Mathematica},
      volume={152},
      number={4},
       pages={876–888},
}

\bib{FL1}{article}{
      author={Friedman, Robert},
      author={Laza, Radu},
       title={Deformations of singular {Fano} and {Calabi}-{Yau} varieties},
        date={2025},
     journal={J. Differential Geom.},
      volume={131},
      number={1},
       pages={65\ndash 131},
}

\bib{HodgeLyubeznik}{article}{
      author={Garc\'{i}a~L\'{o}pez, Ricardo},
      author={Sabbah, Claude},
       title={Hodge–{L}yubeznik numbers},
        date={2025},
     journal={C. R. Math. Acad. Sci. Paris},
      volume={363},
       pages={213\ndash 221},
}

\bib{GNAPGP}{book}{
      author={Guill\'{e}n, F.},
      author={{Navarro Aznar}, V.},
      author={Pascual-Gainza, P.},
      author={Puerta, F.},
       title={Hyperr{\'{e}}solutions cubiques et descente cohomologique},
      series={Lecture Notes in Mathematics},
   publisher={Springer Berlin, Heidelberg},
        date={1988},
      volume={1335},
}

\bib{HartshorneBook}{book}{
      author={Hartshorne, Robin},
       title={{Algebraic} {Geometry}},
      series={Graduate Texts in Mathematics},
   publisher={Springer New York},
        date={2013},
        ISBN={9781475738490},
         url={https://books.google.com/books?id=7z4mBQAAQBAJ},
}

\bib{Hartshorne}{article}{
      author={Hartshorne, Robin},
       title={Ample vector bundles},
        date={1966},
     journal={Publ. Math. Inst. Hautes {\'{E}}tudes Sci.},
      volume={29},
       pages={63\ndash 94},
}

\bib{HartshornePolini}{article}{
      author={Hartshorne, Robin},
      author={Polini, Claudia},
       title={Simple {$\mathscr{D}$}-module components of local cohomology
  modules},
        date={2021},
     journal={J. Algebra},
      volume={571},
       pages={232\ndash 257},
}

\bib{HellusSchenzel}{article}{
      author={Hellus, Michael},
      author={Schenzel, Peter},
       title={On cohomologically complete intersections},
        date={2008},
     journal={J. Algebra},
      volume={320},
      number={10},
       pages={3733\ndash 3748},
}

\bib{DBDeform}{incollection}{
      author={Kov\'{a}cs, S\'{a}ndor},
      author={Schwede, Karl},
       title={Du {B}ois singularities deform},
        date={2016},
   booktitle={Minimal models and extremal rays (kyoto, 2011)},
      volume={70},
   publisher={Math. Soc. Japan},
       pages={49\ndash 65},
}

\bib{Lazarsfeld}{book}{
      author={Lazarsfeld, Robert},
       title={Postivity in {Algebraic} {Geometry} {II}},
      series={Ergebnisse der Mathematik und ihrer Grenzgebiete},
   publisher={Springer-Verlag},
     address={Berlin},
        date={2004},
      volume={49},
}

\bib{MOPW}{article}{
      author={Musta\c{t}\u{a}, Mircea},
      author={Olano, Sebasti\'{a}n},
      author={Popa, Mihnea},
      author={Witaszek, Jakub},
       title={The {Du} {Bois} complex of a hypersurface and the minimal
  exponent},
        date={2023},
     journal={Duke Math. J.},
      volume={173},
      number={2},
       pages={1411–\ndash 1436},
}

\bib{MPLocCoh}{article}{
      author={Musta\c{t}\u{a}, Mircea},
      author={Popa, Mihnea},
       title={Hodge filtration on local cohomology, {D}u {B}ois complex and
  local cohomological dimension},
        date={2022},
     journal={Forum Math. Pi},
      volume={10},
       pages={Paper No. e22, 58},
  url={https://doi-org.myaccess.library.utoronto.ca/10.1017/fmp.2022.15},
      review={\MR{4491455}},
}

\bib{OR}{article}{
      author={Olano, {Sebasti\'{a}n}},
      author={Raychaudhury, Debaditya},
       title={On the local cohomology of secant varieties},
        date={2024},
     journal={preprint arXiv:2407.16688},
}

\bib{PPLefschetz}{article}{
      author={Park, Sung~Gi},
      author={Popa, Mihnea},
       title={Hodge symmetry and {L}efschetz theorems for singular varieties},
        date={2025},
     journal={preprint arXiv: 2410.15638v3},
}

\bib{PPFactorial}{article}{
      author={Park, Sung~Gi},
      author={Popa, Mihnea},
       title={Q-factoriality and {H}odge-{D}u {B}ois theory},
        date={2025},
     journal={preprint arXiv: 2508.17748},
}

\bib{PS}{article}{
      author={Popa, Mihnea},
      author={Shen, Wanchun},
       title={Du {Bois} complexes of cones over singular varieties, local
  cohomological dimension, and {$K$}-groups},
        date={2025},
     journal={Rev. Roumaine Math. Pures Appl.},
      volume={70},
      number={1-2},
       pages={133\ndash 155},
}

\bib{RSW}{article}{
      author={Reichelt, Thomas},
      author={Saito, Morihiko},
      author={Walther, Uli},
       title={Topological calculation of local cohomological dimension},
        date={2023},
     journal={J. Singul.},
      volume={26},
       pages={13\ndash 22},
}

\bib{RW}{article}{
      author={Raicu, Claudiu},
      author={Weyman, Jerzy},
       title={Local cohomology with support in ideals of symmetric minors and
  {P}faffians},
        date={2016},
        ISSN={0024-6107,1469-7750},
     journal={J. Lond. Math. Soc.},
      volume={94},
      number={3},
       pages={709\ndash 725},
         url={https://doi.org/10.1112/jlms/jdw056},
      review={\MR{3614925}},
}

\bib{SaitoMHC}{article}{
      author={Saito, Morihiko},
       title={Mixed {H}odge complexes on algebraic varieties},
        date={2000},
     journal={Math. Ann.},
      volume={316},
      number={2},
       pages={283\ndash 331},
}

\bib{SaitoOnTheHodge}{article}{
      author={Saito, Morihiko},
       title={On the {Hodge} filtration of {Hodge} modules},
        date={2009},
     journal={Mosc. Math. J.},
      volume={9},
      number={1},
       pages={161\ndash 191},
}

\bib{SaitoMHM}{article}{
      author={Saito, Morihiko},
       title={Mixed {H}odge modules},
        date={1990},
        ISSN={0034-5318},
     journal={Publ. Res. Inst. Math. Sci.},
      volume={26},
      number={2},
       pages={221\ndash 333},
         url={https://doi.org/10.2977/prims/1195171082},
      review={\MR{1047415}},
}

\bib{SaitoFormalism}{article}{
      author={Saito, Morihiko},
       title={On the formalism of mixed sheaves},
        date={1991},
     journal={RIMS-preprint},
      number={784},
}

\bib{SaitoCycle}{article}{
      author={Saito, Morihiko},
       title={Refined {C}ycle {M}aps},
        date={2002},
     journal={Advanced Studies in Pure Mathmatics},
      volume={36},
       pages={115\ndash 143},
        note={Algebraic Geometry 2000, Azumino},
}

\bib{SchnellVanishing}{article}{
      author={Schnell, Christian},
       title={On {S}aito's vanishing theorem},
        date={2016},
     journal={Math. Res. Lett.},
      volume={23},
      number={2},
       pages={499\ndash 527},
}

\bib{Tenie}{article}{
      author={Tenie, Anda},
       title={Global smoothing of singular {F}ano and {C}alabi-{Y}au
  varieties},
        date={2026},
     journal={preprint arXiv:2602.07615},
}

\end{biblist}
\end{bibdiv}
\bibliographystyle{abbrv}

\end{document}